\documentclass[10pt]{elsarticle}
\usepackage{amsfonts,enumerate,amsmath,amssymb,mathrsfs,amsbsy}
\usepackage[colorlinks=true]{hyperref}
\usepackage{graphicx}
\usepackage{subfigure}
\usepackage{caption}
\biboptions{numbers,sort&compress}
\usepackage{xcolor}
\usepackage[a4paper,left=30mm,width=150mm,top=25mm,bottom=25mm]{geometry}
\def\qed{\strut\hfill $\Box$}
\newtheorem{thm}{Theorem}[section]

\newtheorem{lem}[thm]{Lemma}

\def\para#1{\vskip .4\baselineskip\noindent{\bf #1}}
\numberwithin{equation}{section}
\allowdisplaybreaks[4]
\begin{document}
	\begin{frontmatter}	
	\title{Averaging principle of stochastic Burgers equation driven by L\'{e}vy processes}
	
	\author[mymainaddress]{Hongge Yue}
	\ead{yuehongge803@163.com}

	\author[mymainaddress,mysecondaryaddress]{Yong Xu}
	\cortext[mycorrespondingauthor]{Corresponding author}
	\ead{hsux3@nwpu.edu.cn}
	
	\author[mymainaddress]{Ruifang Wang}
	\ead{wangruifang_0714@163.com}
	
	\author[mymainaddress]{Zhe Jiao\corref{mycorrespondingauthor}}
	\ead{zjiao@nwpu.edu.cn}
		
	\address[mymainaddress]{School of Mathematics and Stochastics, Northwestern Polytechnical University, Xi'an, 710072, China}

	\address[mysecondaryaddress]{MIIT Key Laboratory of Dynamics and Control of Complex Systems, Northwestern Polytechnical University, Xi'an, 710072, China}

	\begin{abstract}
We are concerned about the averaging principle for the stochastic Burgers equation with slow-fast time scale. This slow-fast system is driven by L\'{e}vy processes. Under some appropriate conditions, we show that the slow component of this system strongly converges to a limit, which is characterized by the solution of stochastic Burgers equation whose coefficients are averaged with respect to the stationary measure of the fast-varying jump-diffusion. To illustrate our theoretical result, we provide some numerical simulations. 
						
	\end{abstract}		
\end{frontmatter}

\section{Introduction }
%

The study of the averaging principle for stochastic partial differential equations with slow-fast time scale has attracted many researchers’ attention (see e.g. \cite{wang2012average, brzezniak2014strong,  2021averagingpei, cerrai2017averaging, fu2017weak, 2018Averagingpei} and the references therein). 
In 2009,  Cerrai  \cite{cerrai2009khasminskii} prove that an averaging principle holds for a general class of stochastic reaction-diffusion systems in any space dimension, which show that the classical Khasminskii approach for systems with a finite number of degrees of freedom can be extended to infinite-dimensional systems.  
Afterwards, an averaging principle for the complex Ginzburg-Landau equations, perturbed by a mixing random force on long time intervals, was established in  \cite{averaging2021gao}.  
Recently, the authors in \cite{ gao2021Averaging}  proved the slow component of the stochastic 2D Navier–Stokes equation converges to the solution of the corresponding averaged equation for any given initial value in the separable real Hilbert space, where the solution is a weak solution.

In this paper, we consider the initial-boundary value problem for the following one dimensional stochastic Burgers equation with slow-fast time scale
\begin{equation}\label{orginal0}
\left\{ 
	\begin{array}{l}
	dX_{t}^{\varepsilon}( \xi ) =\big[\nu\Delta X_{t}^{\varepsilon}( \xi ) +\frac{1}{2}\frac{\partial}{\partial \xi}( X_{t}^{\varepsilon}( \xi ) ) ^2+f_1( X_{t}^{\varepsilon}( \xi ) ,Y_{t}^{\varepsilon}( \xi ) ) \big] dt\\
	\qquad\qquad\quad+dW_{t}^{Q_1}(\xi) +\int_{|z|<1}{h_1( X_{t}^{\varepsilon}( \xi ) ,z )}\tilde{N}_1( dz, dt ) ,\\
	dY_{t}^{\varepsilon}( \xi ) =\frac{1}{\varepsilon}\big[ c\Delta Y_{t}^{\varepsilon}( \xi) +f_2( X_{t}^{\varepsilon}( \xi ) ,Y_{t}^{\varepsilon}( \xi ) ) \big]dt \\
	\qquad\qquad\quad+\frac{1}{\sqrt{\varepsilon}}dW_{t}^{Q_2}(\xi) +\int_{|z|<1}{h_2\big( X_{t}^{\varepsilon}( \xi ) ,Y_{t}^{\varepsilon}( \xi ),z\big)}\tilde{N}_{2}^{\varepsilon}( dz, dt),\\
	X_{0}^{\varepsilon}( \xi) = x ,Y_{0}^{\varepsilon}( \xi ) = y, \quad X_{t}^{\varepsilon}( 0 ) =X_{t}^{\varepsilon}( 1 ) =Y_{t}^{\varepsilon}( 0 ) =Y_{t}^{\varepsilon}( 1) =0,
	\end{array} \right. 
\end{equation}
for $\xi \in [ 0,1 ]$ and $t\in [ 0,T ]$,  $T<\infty$. Here, $\Delta$ is the Laplacian operator, $f_i$ and $h_i$, $i=1, 2$, are nonlinear coefficients. 
The independent Poisson random measures $\tilde{N}_1(\cdot, \cdot)$ and $\tilde{N}_2^\varepsilon(\cdot, \cdot)$, are given by $\tilde{N}_1(dz, dt)=N_1(dz, dt)-\mu_1(dz) dt$, and
$\tilde{N}_2^\varepsilon(dz, dt)=N_2(dz, dt)-\frac{1}{\varepsilon}\mu_2(dz) dt$ respectively,
where $N_i(dz, dt)$, $i=1,2$, is associated Poisson measure, and $\mu_i$, $i=1,2$, the L\'evy measure satisfying $\int_{\mathbb{R} \setminus  0}(1\wedge z^2)\mu_i(dz)< \infty$. 
The $Q_i$-Wiener processes $\{W_t^{Q_i}\}_{t\ge 0}$, $i=1, 2$ are mutually independent with values in $L^2([0,1])$, which are also independent of $N_{i}(\cdot, \cdot)$. 
The scaling parameter $\varepsilon>0 $ is used to describe the separation of time scale between the slow variable $X_t^\varepsilon $ and the fast variable $Y_t^\varepsilon $. $\nu>0$ is the kinematic viscosity and $c\geq 0$ is the diffusion coefficient. For the sake of simplicity, we set the coefficient $\nu$  to be 1. 

This model (\ref{orginal0}) describes Burgers turbulence in the presence of random forces. If the non-Gaussian white noise  in system (\ref{orginal0}) is absent, the authors in   \cite{dong2018averaging} proved the convergence of the slow component, both in the strong sense and in the weak sense. The purpose of this paper is to deal with the case of non-Gaussian white noise. We prove that for any $t\in[0, T]$, as $\varepsilon$ goes to zero, the slow component $X_t^\varepsilon $ of the system (\ref{orginal0}) strongly converges to $\bar{X}_t$ which is the solution of the corresponding average equation. Our proof is based on the Khasminskii argument proposed in \cite{khas1968averaging}. However, L\'evy noise and nonlinear external force, these two terms bring about essential difficulty in proving that the slow component is equicontinuous in some sense. Therefore, we have to find some new ideas to deal with this problem, and obtain a new high-order estimate (see in Lemma \ref{lem3.33}) of the slow component in order to establish the averaging principle for the system (\ref{orginal0}).

The outline of this paper is as follows. In Section \ref{preliminary}, we present some notations, give some precise conditions for the slow-fast system. At the end of this section, we give the statement of our main result.  Section \ref{estimates} gives some prior estimates, which play significant role in proving our main result. In the whole Section \ref{main}, we devote ourselves to the proof of our main result about the averaging principle for the system (\ref{orginal0}). Section \ref{numerical} presents some numerical examples illustrating the averaging principle for the system (\ref{orginal0}).
Throughout the paper, $C, C_i$ are as generic constants whose values may change from line to line.

\section{Preliminaries} \label{preliminary}
The quadruple $(\Omega, \mathcal{F}, \{\mathcal{F}_t\}_{t \geq 0},\mathbb{P})$ is a given stochastic basis satisfying the usual hypotheses in this paper. 
$\mathbb{E}(\cdot)$ stands for expectation with respect to the probability measure $\mathbb{P}$. 
Let $D$ be a bounded domain in a Euclidean space. $L^2(D)$ is the space of square integrable real-valued functions on $D$. The norm in this space is denoted by $\|\cdot\|$. If $k\geq 0$ is an integer, we define $H^{k}(D)$ to consist of all functions in $L^2$ whose differentials belong to $L^2$ up to the order $k$. For real $s\geq 0$, we can define $H^{s}(D)$ by the interpolation $[L^2(D), H^k(D)]_{\theta}$, $k\leq s$ and $s=\theta k$. $H_0^1(D)$ is the subspace of $H^1(D)$ consisting of the functions vanishing on the boundary. Denote by $H^{-1}(D)$ the dual space to $H_0^1(D)$. Throughout this paper, we take $D=(0, 1)$.  
 
Define the bilinear operator $B: L^2\times H_0^1\rightarrow H^{-1}$ by $B(x,y)=x(\xi)\frac{\partial}{\partial {\xi}}y(\xi)$, and the trilinear operator $b: L^2\times H_0^1\times L^2\rightarrow  \mathbb{R}$ by $b( x, y, z ) =\int_0^1{x( \xi )}\frac{\partial y(\xi)}{\partial {\xi}} z( \xi ) d\xi$. Set $B(x)=B(x,y)$ if $x=y$. 

The Laplacian operator is given by
$AX:=\Delta X=\frac{\partial ^2}{\partial \xi ^2}X$ in which $X$ belongs to the domain $\mathscr{D}( A ):=H^2\cap H_{0}^{1}$.
Then from \cite{Guiseppe1994Stochastic} we know that $A$ is the infinitesimal generator of a $C_0$-contraction semigroup $e^{tA}$, $t\geq 0$, which has a regularizing effect, that is, for any $s_1\le s_2 $
\begin{eqnarray}\label{qianru}
\| e^{tA}X \|_{H^{s_2}}\le C\left( 1+t^{\frac{( s_1-s_2 )}{2}} \right) \| X \|_{H^{s_1}}, \quad X\in H^{s_1}.
\end{eqnarray}
The eigenfunctions of $-A$ is given by $e_k =\sqrt{2}\sin(k\pi \xi)$, $k\in\mathbb{N}^+$, $\xi\in[0, 1]$, with the corresponding eigenvalues $\lambda_k=k^2\pi $. 
For any $\alpha\in \mathbb{R}, (-A)^{\alpha}$ is the power of the operator $-A$, and $|\cdot|_\alpha$ is the norm of $ \mathscr{D}\big((-A)^{\frac{\alpha}{2}}\big)$ which is equivalent to the norm of $\mathcal{H}^\alpha$. 

The $Q_1$-Wiener processes $W_t^{Q_1}$ can be given by
\begin{eqnarray}\label{decomposition}
	W_{t}^{Q_1}=\sum_{k=1}^{\infty}{\sqrt{\alpha_k}}\beta _{t}^{k}e_k, t\ge 0,
\end{eqnarray}
where ${\alpha_k}\ge 0$, satisfying $\sum_{k=1}^{\infty}{\alpha_k}<+\infty$, and $\{\beta_t ^k\}_{k\in \mathbb{N}}$ is a sequence of mutually independent standard Brownian motions. And we also assume that $W_t^{Q_2}$ also has a similar decomposition as in (\ref{decomposition}). 

To formulate our main result, we introduce the following assumptions.
\begin{itemize}
	\item[(A1)] There exist two positive constants $L_{f_1}, L_{f_2}$, such that for any $x_1, x_2, y_1,y_2\in L^2$,
	\begin{eqnarray*}
		&&\|f_1(x_1,y_1)\|\leq L_{f_1}(1+\|x_1\|+\|y_1\|),\\
		&&\|f_2(x_1,y_1)\|\leq L_{f_2}(1+\|x_1\|+\|y_1\|),\\	
		&&\|f_1(x_1,y_1)-f_1(x_2,y_2)\|\leq L_{f_1}(\|x_1-x_2\|+\|y_1-y_2\|),\\
		&&\|f_2(x_1,y_1)-f_2(x_2,y_2)\|\leq L_{f_2}(\|x_1-x_2\|+\|y_1-y_2\|).		
	\end{eqnarray*}
\end{itemize}

\begin{itemize}
	\item[(A2)] There exist constants $L_{h_1},L_{h_2}>0 $,  such that for any $\gamma \ge1$, $ \alpha \in \big[1,  \frac{3}{2} \big)$, $x_1, x_2, y_1, y_2\in L^2$, 
	\begin{eqnarray*}
		&&\int_{|z|<1}\|h_1(0,z)\|^\gamma \mu_1(dz) < \infty,~~ \int_{|z|<1}\|h_2(0,0,z)\|^\gamma \mu_2(dz)< \infty,\\
	&&
		\int_{|z|<1}{\lVert h_1( x_1,z ) -h_1( x_2,z) \rVert}^{\gamma }\mu_1( dz ) \leq L_{h_1}\lVert x_1-x_2 \rVert ^{\gamma },
		\\&&
		\int_{|z|<1}{\lVert h_2( x_1, y_1, z ) -h_2( x_2,y_2, z ) \rVert}^{\gamma }\mu_2( dz ) \leq L_{h_2}(\lVert x_1-x_2 \rVert ^{\gamma }+\lVert y_1-y_2 \rVert^{\gamma }),
		\\&&
		\int_{|z|<1} | h_1( x_1,  z ) |_\alpha^{\gamma } \mu_1( dz ) \leq L_{h_1}\left(1+| x_1 |_\alpha ^{\gamma }\right).
	\end{eqnarray*}
\end{itemize}

\begin{itemize}
	\item[(A3)] Let $\eta:=2\lambda _1-L_{f_2}-L_{h_2}>0$. 
\end{itemize} 
\begin{itemize}
\item[(A4)] There exist constants $ \alpha \in \big[1,  \frac{3}{2} \big)$, $\beta\in(0, \infty)$, $\rho\in(2,\infty)$  and $\frac{\beta(\rho-2)}{\rho}<1$  such that  
$$\sum_{k=1}^{\infty}\frac{\alpha_k^{\rho }}{\lambda_k^{\beta}} <+\infty.$$
\end{itemize}

Under the conditions (A1)-(A2), it deduces from \cite{dong2007one} that if the initial data $(x,y)\in H^{\alpha}\times L^2 $, there exist solutions $X_t^{\varepsilon }\in \mathscr{D}( A )$ and $Y_t^{\varepsilon }\in L^2\big( [ 0,T ] ; H_{0}^{1}\big) \cap C\big( [ 0,T ]; L^2  \big)$ satisfying
		\begin{eqnarray}\label{solution}
		\left\{\begin{array}{l}
		X_{t}^{\varepsilon }=e^{tA} x +\int_0^t{e^{( t-s ) A}}B( X_{s}^{\varepsilon } ) ds+\int_0^t{e^{( t-s ) A}}f_1 ( X_{s}^{\varepsilon },Y_{s}^{\varepsilon }  ) ds+\int_0^t{e^{( t-s ) A}}dW_{s}^{Q_1},\\
		\qquad\quad+\int_0^t{\int_{|z|<1}{e^{( t-s  ) A}h_1 ( X_{s}^{\varepsilon },z )}}\tilde{N}_1(ds, dz ),\\
		Y_{t}^{\varepsilon }=e^{\frac{t}{\varepsilon }A} cy +\frac{1}{\varepsilon }\int_0^t{e^{\frac{( t-s ) }{\varepsilon }A}f_2( X_{s}^{\varepsilon },Y_{s}^{\varepsilon } )}ds+\frac{1}{\sqrt{\varepsilon }}\int_0^t{e^{\frac{( t-s ) }{\varepsilon }A}}dW_{s}^{Q_2},\\
		\qquad\quad+\int_0^t{\int_{|z|<1}{e^{\frac{( t-s ) }{\varepsilon }A}h_2( X_{s}^{\varepsilon },Y_{s}^{\varepsilon },z)}}\tilde{N}_2^{\varepsilon }(ds, dz).
		\end{array}
		\right.
		\end{eqnarray} 
For any fixed $x, y\in L^2$,  we consider the following frozen equation associated with the fast component    
\begin{eqnarray}\label{frozeneq}
\left\{\begin{array}{l} 	
d Y_t=[cAY_t +f_2( x,Y_t)]dt +d{W}_{t}^{Q_2}+\int_{|z|<1}{h_2( x,Y_t ,z)}\tilde{N}_2( t,dz),~ Y_0( \xi ) = y ,\\
Y_t( 0 ) =Y_t( 1) =0,~t\in [ 0,T ].
\end{array}
\right.
\end{eqnarray}
The assumptions above imply from  Chapter 16  in \cite{2007Stochastic}  that there exists a unique invariant measure $\mu^{x}$ for (\ref{frozeneq}).

Now we give the statement of our main result.
\begin{thm}\label{thm2.3}
	Under the assumption {\rm (A1)-(A4)},  the slow component $X^{\varepsilon}_t$ of the stochastic Burgers system with the initial data $ (x, y) \in H^\alpha\times L^2$ satisfy the averaging principle
	\[
		\lim_{\varepsilon  \rightarrow 0}\mathbb{E}\Big[\sup\limits_{0\leq t\leq T}\big\|X_t^\varepsilon -\bar{X}_t\big\|^{p}\Big]=0,\quad p\ge 2
	\]
where $\bar{X}_t$ is the solution of the averaged equation
\begin{eqnarray*}
\left\{\begin{array}{l}
d\bar{X}_t=\Delta \bar{X}_tdt+\frac{1}{2}\frac{\partial}{\partial \xi}( \bar{X}_t) ^2dt+\bar{f}_1( \bar{X}_t ) dt+dW_{t}^{Q_1}+\int_{|z|<1}{h_1}( \bar{X}_t , z) \tilde{N}_1( dt, dz) ,\\
\bar{X}_0= x ,
\end{array}
\right.
\end{eqnarray*}
with $\bar{f}_1(x)=\int_{L^2}f_1(x,y)\mu^x(dy)$.
\end{thm}

\section{Prior estimates} \label{estimates}

Before giving the proof of our main result, we need some prior estimates.

\begin{lem}\label{lem3.1}
Assume that the conditions {\rm (A1)-(A2)} are satisfied, then there exist a  positive constant $C_{q,T}$ such that for any $ x \in {H}^\alpha$, $ y \in L^2$, $q \ge 1$, $T>0$ and $  \varepsilon\in(0,1)$
\begin{eqnarray}\label{boundx11}
			\mathbb{E}\left[\underset{0\le t\le T}{\sup} \lVert X_{t}^{\varepsilon } \rVert ^{2q}\right]+\mathbb{E}\left[ \int_0^T\|X_t^\epsilon\|^{2q-2}|X_t^\epsilon |_1^2\right]dt&\le& C_{q,T}( 1+\lVert  x  \rVert ^{2q}+\lVert  y  \rVert ^{2q}),
\end{eqnarray}
	and
	\begin{eqnarray}\label{boundy1}
			\sup \limits_{0\leq t \leq T}\mathbb{E}\| Y_t^\varepsilon \|^{2q}&\le& {C_{q,T}}(1+\| x \|^{2q}+\| y \|^{2q}),
	\end{eqnarray}
where $C_{q,T}$ is independent of $\varepsilon $.
\end{lem}
\para{Proof:} 
According to It\^o's formula, (A1) -(A2) and the Poincaré inequality in \cite{dong2007one}, we have
 \begin{eqnarray}\label{boundy}
\frac{d}{dt}\mathbb{E}\lVert Y_{t}^{\varepsilon } \rVert ^{2q}&=&\frac{2qc}{\varepsilon}\mathbb{E}\left[ \lVert Y_{t}^{\varepsilon } \rVert ^{2q-2}\langle AY_{t}^{\varepsilon },Y_{t}^{\varepsilon }\rangle \right] +\frac{2q}{\varepsilon }\mathbb{E}\left[ \lVert Y_{t}^{\varepsilon } \rVert ^{2q-2}\langle  f_2( X_{t}^{\varepsilon },Y_{t}^{\varepsilon } ) ,Y_{t}^{\varepsilon } \rangle \right] \notag
\\&&
+\frac{q}{\varepsilon }\mathbb{E}\left[\lVert Y_{t}^{\varepsilon } \rVert ^{2q-2}\mathrm{Tr}Q_2\right]+\frac{2q( q-1 )}{\varepsilon }\mathbb{E}\left[\lVert Y_{t}^{\varepsilon } \rVert ^{2q-2}\mathrm{Tr}Q_2\right]\notag
\\&&
+\frac{1}{\varepsilon }\mathbb{E}\left[\int_{|z|<1}{\big( \lVert Y_{t}^{\varepsilon }+h_2( X_{t}^{\varepsilon },Y_{t}^{\varepsilon },z) \rVert ^{2q}-\lVert Y_{t}^{\varepsilon } \rVert ^{2q}\big)}\mu_2( dz )\right] \notag
\\&&
-\frac{2q}{\varepsilon }\mathbb{E}\left[\int_{|z|<1}{\lVert Y_{t}^{\varepsilon } \rVert ^{2q-2}\langle  h_2( X_{t}^{\varepsilon },Y_{t}^{\varepsilon },z ) ,Y_{t}^{\varepsilon } \rangle }\mu_2( dz )\right] \notag\\
&
\le&- \frac{2qc\lambda _1}{\varepsilon }\mathbb{E} \lVert Y_{t}^{\varepsilon } \rVert ^{2q} +\frac{2q}{\varepsilon }\mathbb{E}\big[ \lVert Y_{t}^{\varepsilon } \rVert ^{2q-2}\langle  f_2( X_{t}^{\varepsilon },Y_{t}^{\varepsilon } ) ,Y_{t}^{\varepsilon } \rangle  \big] 
\\&&
+\frac{q}{\varepsilon }\mathbb{E}\left[ \lVert Y_{t}^{\varepsilon } \rVert ^{2q-2}\mathrm{Tr}Q_2\right]+\frac{2q( q-1 )}{\varepsilon }\mathbb{E}\left[\lVert Y_{t}^{\varepsilon } \rVert ^{2q-2}\mathrm{Tr}Q_2\right]\notag
\\&&
+\frac{1}{\varepsilon }\mathbb{E}\left[\sum_{i=2}^{2q}{C_{2q}^{i}}\int_{|z|<1}{\lVert Y_{t}^{\varepsilon } \rVert ^{2q-i}\lVert h_2( X_{t}^{\varepsilon },Y_{t}^{\varepsilon },z )  \rVert ^i}\mu_2( dz)\right] \notag
\end{eqnarray}
From (\ref{boundy}), using condition (A3), and the Young's inequality, we deduce that there exist a positive constant $\gamma'$ such that
\begin{eqnarray} \label{boundy2}
\frac{d}{dt}\mathbb{E}\lVert Y_{t}^{\varepsilon } \rVert ^{2q}\le -\frac{q\gamma'}{\varepsilon }\mathbb{E}\lVert Y_{t}^{\varepsilon } \rVert ^{2q}+\frac{C_q}{\varepsilon }\mathbb{E}\lVert X_{t}^{\varepsilon } \rVert ^{2q}+\frac{C_q}{\varepsilon }.
\end{eqnarray}
Applying the comparison theorem, we get
\begin{eqnarray} \label{boundy3}
\mathbb{E}\lVert Y_{t}^{\varepsilon } \rVert ^{2q}\le \lVert  y  \rVert ^{2q}e^{-\frac{q\gamma'}{\varepsilon }t}+\frac{C_q}{\varepsilon }\int_0^t{e^{-\frac{q\gamma}{\varepsilon }( t-s)}( 1+\mathbb{E}\lVert X_{s}^{\varepsilon } \rVert ^{2q} ) ds}.
\end{eqnarray}
For $X_t^\varepsilon$, by It\^o's formula, we have
\begin{eqnarray} \label{boundy2}
\mathbb{E}\lVert X_{t}^{\varepsilon } \rVert ^{2q}:=\|x\|^{2q}+\sum  \limits_{i=1}^9 \varXi_i(t),
\end{eqnarray}
where 
\begin{eqnarray*}
 \varXi_1(t)&=&2q\mathbb{E}\left[ \int_0^t \left[ \lVert X_{s}^{\varepsilon } \rVert ^{2q-2}\langle AX_{s}^{\varepsilon },X_{s}^{\varepsilon } \rangle \right]ds\right],\\
 \varXi_2(t)&=&2q\mathbb{E}\left[\int_0^t  \left[ \lVert X_{s}^{\varepsilon } \rVert ^{2q-2}\langle  B(X_{s}^{\varepsilon }),X_{s}^{\varepsilon } \rangle \right]ds\right],\\
 \varXi_3(t)&=&2q\mathbb{E}\left[\int_0^t \left[ \lVert Y_{s}^{\varepsilon } \rVert ^{2q-2}\langle  f_1( X_{s}^{\varepsilon },Y_{s}^{\varepsilon } ) ,X_{s}^{\varepsilon } \rangle  \right] ds\right],\\
 \varXi_4(t)&=&q\mathbb{E}\left[\int_0^t \lVert X_{s}^{\varepsilon } \rVert ^{2q-2}\mathrm{Tr}Q_1ds\right],\\
 \varXi_5(t)&=&2q( q-1 ) \mathbb{E}\left[\int_0^t\lVert X_{s}^{\varepsilon } \rVert ^{2q-2}\mathrm{Tr}Q_1ds\right],\\
 \varXi_6(t)&=&2q\mathbb{E}\left[\int_0^t \lVert X_{s}^{\varepsilon } \rVert ^{2q-2}  \cdot \lVert X_{s}^{\varepsilon } \rVert dW_s^{Q_1}\right],\\
 \varXi_7(t)&=&\mathbb{E}\left[\int_0^t\int_{|z|<1}{\big(\lVert X_{s}^{\varepsilon }+h_1( X_{s}^{\varepsilon },Y_{s}^{\varepsilon },z ) \rVert ^{2q}-\lVert X_{s}^{\varepsilon } \rVert ^{2q}\big)} \tilde{N}_1( ds, dz)\right],\\
 \varXi_8(t)&=&\mathbb{E}\left[\int_0^t\int_{|z|<1}{\big(\lVert X_{s}^{\varepsilon }+h_1( X_{s}^{\varepsilon },Y_{s}^{\varepsilon },z ) \rVert ^{2q}-\lVert X_{s}^{\varepsilon } \rVert ^{2q}\big)}\mu_1( dz )ds\right],\\
 \varXi_9(t)&=&-2q\mathbb{E}\left[\int_0^t\int_{|z|<1}{\lVert X_{s}^{\varepsilon } \rVert ^{2q-2}\langle  h_1( X_{s}^{\varepsilon },Y_{s}^{\varepsilon },z) ,X_{s}^{\varepsilon } \rangle }\mu_1( dz )ds\right].
\end{eqnarray*}
 For any $x,y\in  {H}_0^1$,  $b(x,y,y)=0$( see e.g. \cite{dong2007one}, Proposition 2.3.),  we get
\begin{eqnarray*}
\mathbb{E}\left[\int_0^t  \left[ \lVert X_{s}^{\varepsilon } \rVert ^{2q-2}\langle  B(X_{s}^{\varepsilon }),X_{s}^{\varepsilon } \rangle \right]ds\right]
&\leq &\sup_{0\leq t\leq T} \lVert X_{t}^{\varepsilon } \rVert ^{2q-2} \left[\int_0^t  \langle  B(X_{s}^{\varepsilon }),X_{s}^{\varepsilon } \rangle ds\right]\\
&\leq &C_T \sup_{0\leq t\leq T} \lVert X_{t}^{\varepsilon } \rVert ^{2q-2} \left[\int_0^1 \langle  B(X_{sT}^{\varepsilon }),X_{sT}^{\varepsilon } \rangle ds\right]\\
&\leq &C_T \sup_{0\leq t\leq T} \lVert X_{t}^{\varepsilon } \rVert ^{2q-2} b(X_{sT}^{\varepsilon } ,X_{sT}^{\varepsilon } ,X_{sT}^{\varepsilon } )\\
&=&0.
\end{eqnarray*}
Note that $\left\langle AX, X\right\rangle  =-\| \nabla X\|^2\le-\beta|X|_1^2+\gamma\|X\| ^2 $, $\beta>0, \gamma>0$ (see e.g. \cite{dong2018averaging}),   by the linear growth conditions in (A1)  and Young's inequality, we have
\begin{eqnarray} \label{5eq123}
\mathbb{E}\left[\sum_{i=1}^{5}\sup\limits_{0\leq t \leq T}\varXi
_i(t) \right]
&\leq &-2q\int_0^t\|X_s^\epsilon\|^{2q-2}|X_s^\epsilon|_1^2ds\notag\\
&&+C_{q,T}\mathbb{E}\left[ \int_0^T\|X_t^\varepsilon\|^{2q}dt\right]+C_{q,T}\mathbb{E}\left[ \int_0^T\|Y_t^\varepsilon\|^{2q}dt\right]+C_{q,T}.
\end{eqnarray}
On the one hand, with the help of Burkholder-Davis-Gundy inequality, it yields
\begin{eqnarray}\label{5eq2}
\mathbb{E}\left[\sup\limits_{0\leq t \leq T}\varXi
_6(t)\right]&\leq &C_q\mathbb{E}\left[\left(\int_0^T\|X_{t}^\varepsilon\|^{4q-2}\mathrm{Tr}Q_1^2dt\right)^\frac{1}{2}\right]\notag\\
&\leq & C_q \mathbb{E}\left[\left(\sup\limits_{0\leq t \leq T}\|X_{t}^\varepsilon\|^{2q}\int_0^T\|X_{t}^\varepsilon\|^{2q-2}  dt\right)^\frac{1}{2}\right]\\
&\leq & C_{q,T}+C_{q,T}\int_0^T\mathbb{E}\|X_{t}^\varepsilon\|^{2q}dt+\frac{1}{3}\mathbb{E}\sup\limits_{0\leq t \leq T}\|X_{t}^\varepsilon\|^{2q}\notag.
\end{eqnarray}
On the other hand, by Burkholder-Davis-Gundy type inequality for stochastic integral with respect to Poisson compensated martingale measures and Young inequality, we have
\begin{eqnarray}\label{5eq3}
\mathbb{E}\left[\sup\limits_{0\leq t \leq T}\varXi
_7(t)\right]&=&\mathbb{E} \left[ \sup_{0\leq t\leq T}\int^{t}_{0} \int_{{|z|<1}} \left[ \|X^{\varepsilon }_{s}+h_{1}\left( X^{\varepsilon }_{s},z\right)  \|^{2q}-\|X^{\varepsilon }_{s}\|^{2q}\right]  \tilde{N_{1}} (dz,ds)\right] \notag \\
&\leq& C\mathbb{E} \left[\left( \int^{T}_{0} \int_{{|z|<1}} \left(\| X^{\varepsilon }_{t}+h_{1}\left( X^{\varepsilon }_{t},z\right)  \|^{2q}-\|X^{\varepsilon }_{t}\|^{2q}\right)^{2}  \mu_1(dz)dt\right)^{\frac{1}{2} } \right] \notag\\
&\leq &C\mathbb{E} \left[ \int^{T}_{0} \int_{{|z|<1}} (\| X^{\varepsilon }_{t}\| ^{4q-2}\| h_{1}\left( X^{\varepsilon }_{t},z\right)  \|^{2} +\| h_{1}\left( X^{\varepsilon }_{t},z\right)  \|^{4q} )\mu_1(dz)dt)\right]^{\frac{1}{2} }  \notag\\
&\leq& C\mathbb{E} \left[ \left( \sup_{0\leq t\leq T} \| X^{\varepsilon }_{t}\|^{2q}   \int^{T}_{0} \int_{{|z|<1}} \|X^{\varepsilon }_t\|^{2q-2}\| h_{1}\left( X^{\varepsilon }_{t},z\right)  \|^{2} \mu_1(dz)dt\right) ^{\frac{1}{2} }\right]\notag\\
&&+C\mathbb{E} \left[ \left(\int^{T}_{0} \int_{{|z|<1}} \| h_{1}\left( X^{\varepsilon }_{t},z\right)  \|^{4q} \mu_1(dz)dt)\right)^{\frac{1}{2} } \right] \\
&\leq&\frac{1}{3} \mathbb{E} \left[\sup_{0\leq t\leq T} \| X^{\varepsilon }_{t}\|^{2q} \right] +C_{q,T}\mathbb{E} \left[\int^{T}_{0} \int_{{|z|<1}} \| h_{1}\left( X^{\varepsilon }_{t},z\right) \| ^{2q} \mu_{1} \left( dz\right)  dt\right]\notag\\
&&+C_{q,T}\mathbb{E} \left[\int^{T}_{0} \int_{{|z|<1}} \|X^{\varepsilon }_t\|^{2q-2}\| h_{1}\left( X^{\varepsilon }_{t},z\right) \| ^2 \mu_{1} \left( dz\right)  dt\right]\notag \\
&\leq&\frac{1}{3}\mathbb{E} \left[ \sup_{0\leq t\leq T} \| X^{\varepsilon }_{t}\|^{2q} \right]  +C_{q,T}\int^{T}_{0} \mathbb{E} \| X^{\varepsilon }_{t}\|^{2q} dt+C_{q,T}\notag.
\end{eqnarray}
Next,  by Binomial theorem, (A2) and Young's inequality, we have
\begin{eqnarray} \label{5eq1}
\mathbb{E}\left[\sum_{i=8}^9\sup\limits_{0\leq t \leq T}\varXi
_i(t)\right]&=&\mathbb{E}\left[\sup\limits_{0\leq t \leq T}\int_0^t\int_{|z|<1}\sum \limits_{i=2}^{2q}C_{2q}^i\|X_{s}^\varepsilon\|^{2q-i}\|h_1(X_{s}^\varepsilon,z)\|^i \mu_{1}(dz)ds\right]\notag\\
&\leq& C_q\sum \limits_{i=2}^{2q}C_{2q}^i\mathbb{E}\left[\sup\limits_{0\leq t \leq T}\int_0^t\|X_{s}^\varepsilon\|^{2q-i}(1+\|X_{s}^\varepsilon\|^{i})ds\right]\\
&\leq& C_{q,T}+C_{q,T}\mathbb{E}\int_0^T\|X_{t}^\varepsilon\|^{2q}dt\notag.
\end{eqnarray}
Combining the estimates  (\ref{5eq123})  (\ref{5eq2}) (\ref{5eq3})  (\ref{5eq1}) and (\ref{boundy3}), it is easy to see that
\begin{eqnarray*}
	&&\mathbb{E}\left[\sup \limits_{0\leq t \leq T}\| X_t^\varepsilon\|^{2q}\right]+2q\mathbb{E}\left[ \int_0^T\|X_t^\epsilon\|^{2q-2}|X_t^\epsilon |_1^2\right]dt
	\\&\leq &C_{q,T}(1+\| x \|^{2q})+C_{q,T}\int_0^{T}\mathbb{E}\left[\sup\limits_{0\leq r\leq s}\|X_r^\varepsilon\|^{2q}\right]ds+C_{q,T}\int_0^{T}\mathbb{E}\left[\sup\limits_{0\leq r\leq s}\|Y_r^\varepsilon\|^{2q}\right]ds\\
	&\leq&C_{q,T}(1+\| x \|^{2q}+\| y \|^{2q})+C_{q,T}\int_0^{T}\mathbb{E}\left[\sup\limits_{0\leq r\leq t}\|X_r^\varepsilon\|^{2q}\right]dt\\
	&&+\frac{C_q}{\varepsilon }\int_0^{T}\int_0^t{e^{-\frac{q\gamma}{\varepsilon }( t-s)}( 1+\mathbb{E}\lVert X_{s}^{\varepsilon } \rVert ^{2q} ) ds} dt\\
	&\leq&C_{q,T}(1+\| x \|^{2q}+\| y \|^{2q})+C_{q,T}\int_0^{T}\mathbb{E}\left[\sup\limits_{0\leq r\leq t}\|X_r^\varepsilon\|^{2q}\right]dt.
\end{eqnarray*}
It then follows from Gronwall's inequality, 
\begin{eqnarray*}
\mathbb{E}\left[\sup \limits_{0\leq t \leq T}\| X_t^\varepsilon\|^{2q}\right]+2q\mathbb{E}\left[ \int_0^T\|X_t^\epsilon\|^{2q-2}|X_t^\epsilon |_1^2\right]dt\leq C_{q,T}(1+\| x \|^{2q}+\| y \|^{2q}),
\end{eqnarray*}	
which also gives
\begin{eqnarray*}
\sup \limits_{0\leq t \leq T} \|Y_t^\varepsilon\|^{2q}\leq C_{q,T}(1+\| x \|^{2q}+\| y \|^{2q}).
\end{eqnarray*}

The proof is completed.    \qed
\begin{lem}\label{lem3.33}
	Assume that the conditions {\rm (A1)-(A4)} are satisfied, then there exist a  positive constant $C_{q,T,\alpha}$ such that for any $ x \in {H}^\alpha$, $ y \in L^2,$  $q>2$, $\alpha \in \big[1,  \frac{3}{2} \big)$ and  $ \varepsilon\in(0,1) $		
\begin{eqnarray}\label{key}
	\mathbb{E}\left[\sup_{0\leq t\leq T }| X_{t}^{\varepsilon }|_{\alpha}^{2q}\right]\leq C_{q,T,\alpha,n}( 1+|  x  |_\alpha^{2q}+\lVert  y  \rVert ^{2q}),
\end{eqnarray}	
 where $C_{q,T,\alpha}$ is independent of $\varepsilon$.
\end{lem}
\para{Proof:}  From Eq.(\ref{solution}), we note $X_{t}^{\varepsilon }=:\sum_{i=1}^5J_i$,
in which
\begin{eqnarray*}
	J_1&=&e^{tA} x,\\
	J_2&=&\int_0^t{e}^{( t-s) A}B( X_{s}^{\varepsilon } ) ds,\\
	J_3&=&\int_0^t{e}^{( t-s) A}f_1( X_{s}^{\varepsilon },Y_{s}^{\varepsilon } ) ds,\\
	J_4&=&\int_0^t{e}^{( t-s) A} dW_{s}^{Q_1},\\
	J_5&=&\int_0^t{\int_{|z|<1}{e^{( t-s) A} h_1(X_{s}^{\varepsilon },z)}}\tilde{N}_1( ds,dz).	
\end{eqnarray*}	

For $J_1$, it is clear that  
\begin{eqnarray}\label{eat1}
| e^{tA} x |_{\alpha}^{2q}\le |  x  |_{\alpha}^{2q}.
\end{eqnarray}

For $J _2$, according to (\ref{qianru}),  Lemma (\ref{lem3.1}) and Sobolev inequalities in ({\cite{temam1983navier}, Lemma 2.1 ), we obtain
\begin{eqnarray*}
\mathbb{E}\underset{0\le t\le T}{\text{sup}}\left| \int_0^t{e}^{\left( t-s \right) A}B\left( X_{s}^{\epsilon} \right) ds \right|_{\alpha}^{2q}
&\le&C\mathbb{E}\left[ \underset{0\le t\le T}{\text{sup}}\int_0^t{\left( 1+\left( t-s \right) ^{\frac{-\alpha _3-\alpha}{2}} \right)}|B(X_s^\epsilon)|_{-\alpha_3} ds \right] ^{2q}\notag\\
&\le&C\mathbb{E}\left[ \underset{0\le t\le T}{\text{sup}}\int_0^t{\left( 1+\left( t-s \right) ^{\frac{-\alpha _3-\alpha}{2}} \right)}|X_s^\epsilon|_{\alpha_1}|X_s^\epsilon|_{\alpha_2+1} ds \right] ^{2q}\notag,
\end{eqnarray*}
where  $\alpha_1+\alpha_2+\alpha_3>\frac{1}{2}, \alpha_i> 0(i=1,2,3)$. Using the  interpolation inequality, we have that
\begin{eqnarray}\label{1}
|X_s^\epsilon|_{\alpha_1}\leq C\|X_s^\varepsilon\|^{\frac{\alpha-\alpha_1}{\alpha}} |X_s^\varepsilon|_{\alpha}^{\frac{\alpha_1}{\alpha}},
\end{eqnarray}
for any $0< \alpha_1< \alpha$, and that
\begin{eqnarray}\label{2}
|X_s^\epsilon|_{\alpha_2+1}\leq C\|X_s^\varepsilon\|^{\frac{\alpha-\alpha_2-1}{\alpha}} |X_s^\varepsilon|_{\alpha}^{\frac{\alpha_2+1}{\alpha}},
\end{eqnarray}
for any $0< \alpha_2+1<\alpha.$ Let $\alpha_1$ and $\alpha_2$ be small enough such that $1+\alpha_1+\alpha_2\in(1, \alpha).$  It follows from (\ref{1}) and (\ref{2}) ,  and according to H\"{o}lder inequality, Young's inequality,  we obtain that 
\begin{eqnarray}\label{second1}
&&\mathbb{E}\underset{0\le t\le T}{\text{sup}}\left| \int_0^t{e}^{\left( t-s \right) A}B\left( X_{s}^{\epsilon} \right) ds \right|_{\alpha}^{2q}\notag\\
&\leq& C\mathbb{E}\left[ \underset{0\le t\le T}{\text{sup}}\int_0^t{\left( 1+\left( t-s \right) ^{\frac{-\alpha _3-\alpha}{2}} \right)}\lVert X_{s}^{\epsilon} \rVert ^{\frac{2\alpha -\alpha _1-\alpha _2-1}{\alpha}}\left| X_{s}^{\epsilon} \right|_{\alpha}^{\frac{\alpha _1+\alpha _2+1}{2}}ds \right] ^{2q}\notag
\\&\leq&
 C\mathbb{E}\left[ \underset{0\le t\le T}{\text{sup}}\int_0^t{\left( 1+\left( t-s \right) ^{\frac{-\alpha _3-\alpha}{2}} \right)}^{\frac{2q}{2q-1}}ds \right] ^{2q-1} \\
 &&\times \mathbb{E}\left[ \underset{0\le t\le T}{\text{sup}}\int_0^t{\lVert X_{s}^{\epsilon} \rVert ^{2q\cdot\frac{2\alpha -\alpha _1-\alpha _2-1}{\alpha}}}\left| X_{s}^{\epsilon} \right|_{\alpha}^{2q\cdot\frac{\alpha _1+\alpha _2+1}{2}}ds \right]\notag
\\&\leq&
 C_{q,T}\left( 1+\int_0^T{s^{\frac{-\alpha _3-\alpha}{2}\cdot\frac{2q}{2q-1}}}ds \right) ^{2q-1}\left( \int_0^T{\mathbb{E}\lVert X_{s}^{\epsilon} \rVert ^{2q\cdot\frac{2\alpha -\alpha _1-\alpha _2-1}{\alpha -\alpha _1-\alpha _2-1}}}ds+\int_0^T{\mathbb{E}\left| X_{s}^{\epsilon} \right|_{\alpha}^{2q}}ds \right)\notag\\
&\le&
 C_{q,T}\left( 1+\int_0^T{\mathbb{E}\left| X_{s}^{\epsilon} \right|_{\alpha}^{2q}}ds \right),\notag
\end{eqnarray}
 let $q$ be large enough such that  ${\frac{\alpha _3+\alpha}{2}\cdot \frac{2q}{2q-1}}<1, \alpha\in(1, \frac{3}{2}).$  Next, let $\alpha=1,$ we have
 \begin{eqnarray}\label{second2}
&&\mathbb{E} \left[ \sup_{0\leq t\leq T} |\int^{t}_{0} {e}^{(t-s)A}  B(X^{\varepsilon }_{s})ds|^{2q}_{1 }\right]  \notag\\
&&\leq \mathbb{E} \left[\left(\sup_{0\leq t\leq T} \int^{t}_{0} \left| {e}^{(t-s)A}  B(X^{\varepsilon }_{s})\right|_{1 }  ds\right)^{2q} \right]\notag\\
&& \leq \mathbb{E}  \left[\left( \sup_{0\leq t\leq T} \int^{t}_{0} \left( 1+\left( t-s\right)^{-\frac{\alpha_3'+1}{2} }  \right)  \left| B(X^{\varepsilon }_{s})\right|_{-\alpha_3'}  ds\right)^{2q} \right]\notag
\\&& \leq C\mathbb{E}  \left[\left( \sup_{0\leq t\leq T} \int^{t}_{0} \left( 1+\left( t-s\right)^{-\frac{\alpha_3'+1}{2} }  \right)  \| X^{\varepsilon }_{s}\| \left| X^{\varepsilon }_{s}\right|_{1}  ds\right)^{2q}  \right]\\
&&\leq C_q\mathbb{E} \left\{ \sup_{0\leq t \leq T}\left [\left(\int^t_{0} \left( 1+\left( t-s\right)^{-\frac{\alpha_3' +1}{2} }  \right)^{\frac{2q}{2q-1} }  ds\right)^{2q-1} \int^t_{0}\| X^{\varepsilon }_{s}\|^{2q} | X^{\varepsilon }_{s}|_1^{2q}ds\right]\right\}\notag\\
&&\leq C_q\left [\left(\int^T_{0} \left( 1+\left( t-s\right)^{-\frac{\alpha_3'+1}{2} }  \right)^{\frac{2q}{2q-1} }  ds\right)^{2q-1}\mathbb{E}\left[ \sup_{0\leq t\leq T}\| X^{\varepsilon }_{t}\|^{2q}\right] \int^t_{0} \mathbb{E} | X^{\varepsilon }_{s}|_1^{2q}ds\right]\notag\\
&&\leq C_{q, T}(1+\|x\|^{2q}+\|y\|^{2q})\int^{T}_{0} \mathbb{E} | X^{\varepsilon }_{s}|_1^{2q} ds,\notag
\end{eqnarray}
for $ q>2$, where $\alpha_3'\in(\frac{1}{2}, \frac{2}{3}).$

For $J _3$, according to (\ref{qianru}), taking $q$ large enough such that $\alpha\in [0,2-\frac{1}{q})$, it follows from Lemma \ref{lem3.1} that 
\begin{eqnarray}\label{third}
&&\mathbb{E}\Big[\underset{0\le t\le T}{\text{sup}}\big| \int_0^t{e}^{( t-s) A}f_1( X_{s}^{\varepsilon },Y_{s}^{\varepsilon } ) ds \big|_{\alpha}^{2q}\Big]\notag\\
&\le&C\mathbb{E}  \left[\left( \sup_{0\leq t\leq T} \int^{t}_{0} \left( 1+\left( t-s\right)^{-\frac{\alpha }{2} }  \right)  \left( 1+\| X^{\varepsilon }_{s}\| +\| Y^{\varepsilon }_{s}\| \right)  ds\right)^{2q}\right] \notag \\
&\le& C_{q,T}\left( 1+\int^{T}_{0} s^{-\frac{\alpha q}{2q-1} }ds\right)^{2q-1} \mathbb{E}\left[\int^{T}_{0} \left( 1+\| X^{\varepsilon }_{s}\|^{2q} + \| Y^{\varepsilon }_{s}\|^{2q} \right)  ds\right]\\
&\le& C_{q,T}{\big( 1+\lVert  x  \rVert ^{2q}+\lVert  y  \rVert ^{2q} \big)}.\notag
\end{eqnarray}	

For $J _4$, following the same argument as in the proof of Lemma 4.1 in \cite{cerrai2009khasminskii}, according to (A4),  and $\bar{\alpha}>0$ such that for any $\alpha\in[ 0,\bar{\alpha})$,  $q>1$ and  $\varepsilon\in(0,1)$, we have
	\begin{eqnarray}\label{nat}
	\mathbb{E}\Big[\sup_{0\le t\le T}|\int_0^t{{e^{\left( t-s \right) A}}} dW_s^{Q_1}  |_{\alpha}^{2q}\Big]\leq   C_{q,T,\alpha}{\big( 1+\lVert  x  \rVert ^{2q}+\lVert  y  \rVert ^{2q} \big)}, 
	\end{eqnarray}	
	 choose $\theta>0$ such that $2\theta+\frac{2\beta(\rho-2)+{\alpha}(\rho+2)}{2\rho}<1$. 
	 Here, the detailed proof is omitted for the sake of brevity.	 

For $J _5$, 
with aid of the Kunita's first inequality \cite{Applebaum2009Levy}, H\"{o}lder's inequality, (A2), (\ref{qianru})  and (\ref{boundx11}), there holds
 \begin{eqnarray} \label{jump}
&&\mathbb{E} \left[ \sup_{0\leq t\leq T} |\int^{t}_{0} \int_{{|z|<1} } e^{\left( t-s\right)  A}h_{1}\left( X^{\varepsilon }_{s},z\right)  \tilde{N}_{1} \left( ds,dz\right)  |^{2q}_{\alpha }\right]  \notag\\
&\leq&\mathbb{E} \left[ \sup_{0\leq t\leq T}\left(\int^{t}_{0} \int_{{|z|<1} }  |e^{\left( t-s\right)  A}h_{1}\left( X^{\varepsilon }_{s},z\right)    |_\alpha \tilde{N}_{1} \left( ds,dz\right)\right) ^{2q}\right]  \notag\\
&\leq&C\mathbb{E} \left[\left(\int^T_{0} \int_{{|z|<1} }  |e^{\left( t-s\right)  A}h_{1}\left( X^{\varepsilon }_{s}, z\right)    |_\alpha^2 \mu_{1} \left( dz\right)  ds\right)^q\right]\\
&& +C\mathbb{E} \left[\int^{T}_{0} \int_{{|z|<1} }  |e^{\left( t-s\right)  A}h_{1}\left( X^{\varepsilon }_{s}, z\right)  |_\alpha^{2q} \mu_{1} \left( dz\right)  ds \right]  \notag\\
&\leq&C\mathbb{E} \left[\left(\int^T_{0}  (1+| X^{\varepsilon }_{s}| _\alpha^2   ) ds\right)^q\right] +C\mathbb{E} \left[\int^T_{0}(1+| X^{\varepsilon }_{s}| _\alpha^{2q}   ) ds \right]  \notag\\
&\leq&C_{q,T} \left(1+
 \int^{T}_{0} \mathbb{E} |X^{\varepsilon }_{s}  |_\alpha^{2q} ds \right).\notag
\end{eqnarray}

 We conclude the proof by combining (\ref{eat1}),  (\ref{second1}),  (\ref{second2}),  (\ref{third}),  (\ref{nat}), (\ref{jump}) and Gronwall's inequality. \qed
 
\begin{lem}\label{lem3.44}
Assume that the conditions {\rm (A1)-(A4)} are satisfied, then there exist a  positive constant $C_{q,T,\alpha}$ such that for any $ x \in \mathit{H}^\alpha,$ $ y  \in \mathit{L}^2$,  $0\leq  t\leq t+h\leq T$,  $\alpha \in \big[1, \frac{3}{2} \big )$ and  $ \varepsilon\in(0,1)$
\begin{eqnarray*}
\mathbb{E}\left[\big\lVert X_{t+h}^{\varepsilon }-X_{t}^{\varepsilon }\big \rVert ^{2q}\right]\le C_{q,T,\alpha}h^{\alpha q}\big( 1+|  x  |_{\alpha} ^{2q}+\lVert  y   \rVert  ^{2q} \big),
\end{eqnarray*}	
 where $C_{q,T,\alpha}$ is independent of $\varepsilon$.
\end{lem}
\para{Proof:} After simple calculations, we have $X_{t+h}^{\varepsilon }-X_{t}^{\varepsilon } =:\sum_{i=1}^5 \Theta _i$, in which
\begin{eqnarray*}\label{i1}
 \Theta _1&=&( e^{Ah}-I ) X_{t}^{\varepsilon }, \\
  \Theta _2&=&\int_t^{t+h}{e^{( t+h-s ) A}}B( X_{s}^{\varepsilon }) ds,\\
   \Theta _3&=&\int_t^{t+h}{e^{( t+h-s) A}}f_1( X_{s}^{\varepsilon },Y_{s}^{\varepsilon } ) ds, \\
    \Theta _4&=&\int_t^{t+h}{e^{( t+h-s) A}}dW_{s}^{Q_1},\\
     \Theta _5&=&\int_t^{t+h}{\int_{|z|<1}{e^{( t+h-s) A}}}h_1( X_{s}^{\varepsilon },z )\tilde{N}_1(ds, dz ).
\end{eqnarray*}

For $ \Theta _1$, There exist a constant $C_\alpha>0$ such that for any $X\in \mathscr{D}\big((-A)^{\frac{\alpha}{2}}\big),$ $\|e^{Ah}X-X\|\leq C_\alpha h^{\frac{\alpha}{2}}|X|_{\alpha}$ (see e.g. \cite{2012Semigroups}). Then using Lemma {\ref{lem3.33}}, we get
\begin{eqnarray}\label{theta1}
\mathbb{E}\left[\|\Theta _1\|^{2q}\right]&\leq &C_\alpha h^{\alpha q}\mathbb{E}\left[|X_t^{\alpha}|_{\alpha}^{2q}\right]\notag \\
&\leq& C_{q,T,\alpha} h^{\alpha q}(1+| x |_{\alpha}^{2q}+\| y \|^{2q}).
\end{eqnarray}

For $ \Theta _2$, using the contractive property of the semigroup $e^{tA}$, Lemma \ref{lem3.1} and  Lemma \ref{lem3.33} , we obtain
\begin{eqnarray}\label{theta2}
\mathbb{E}\left[\big\lVert \int_t^{t+h}{e^{( t+h-s ) A}}B( X_{s}^{\varepsilon } ) ds \big\rVert ^{2q}\right]&\le& \mathbb{E}\Big[ \int_t^{t+h}\lVert B( X_{s}^{\varepsilon }) \rVert ds \Big] ^{2q}\notag
\\&\le& C_{q,T}\mathbb{E}\left[\left( \int_t^{t+h}{| X_{s}^{\varepsilon } |}_{1}^{2}ds \right)^{2q} \right]\notag
\\
&\le& C_{q,T}h^{2q}\notag.
\end{eqnarray}

For $ \Theta _3$, applying (A1) and Lemma \ref{lem3.1}, we get
\begin{eqnarray}\label{theta3}
\mathbb{E} \| \Theta_{3} \|^{2q} \leq C h^{2q-1}\mathbb{E}\left[ \int^{t+h}_{t} \| f\left( X^{\varepsilon }_{s},Y^{\varepsilon }_{s}\right)  \|^{2q} ds\right]\leq C_{q,T}h^{2q}\left( 1+\| x\|^{2q} +\| y\|^{2q} \right) . 
\end{eqnarray}

For $ \Theta _4$, note that $\Theta _4$ is the centered Gaussian random variable with the variance given by $S_h =\int_o^te^{{h-r}A}Q_1e^{(h-r)A^*}dr $. Then, for any $q > 1,$ we get
\begin{eqnarray}\label{theta4}
\mathbb{E} \| \Theta_{4} \|^{2q} \leq C_{q,T}\left[ \mathrm{Tr} \left( S_{h}\right)  \right]^{q}  =C_q\left( \sum^{\infty }_{k=1} \int^{h}_{0} e^{-2\left( h-r\right)  \lambda_{k} }\alpha_{k} dr\right)^{q}  \leq C_q\left( \sum^{\infty }_{k=1} \alpha_{k} \right)^{q}  h^{q}.
\end{eqnarray}

For $ \Theta _5$,  applying Burkholder-Davis-Gundy inequality  \cite{marinelli2016maximal},  H\"{o}lder's inequality and  (A2), it yields
\begin{eqnarray}\label{theta5}
&&\mathbb{E}\left[\big\lVert \int_t^{t+h}{\int_{|z|<1}{e^{( t+h-s ) A}}}h_1( X_{s}^{\varepsilon },z )\tilde{N}_1(ds, dz )\big\rVert ^{2q}\right]\notag\\
&
\le& C_q\mathbb{E}\left[ \left(\int_t^{t+h}{\int_{|z|<1}{\lVert e^{( t+h-s ) A}h_1( X_{s}^{\varepsilon },z ) \rVert^2}}\mu_1( dz ) ds\right) ^{q}\right]\notag\\
&&+C_q\mathbb{E}\left[ \int_t^{t+h}{\int_{|z|<1}{\lVert e^{( t+h-s ) A}h_1( X_{s}^{\varepsilon },z ) \rVert^{2q}}}\mu_1( dz ) ds\right] 
\\&
\le& C_qh^{q-1}\mathbb{E}\left[\int_t^{t+h}{( 1+\lVert X_{s}^{\varepsilon } \rVert ^{2} )^q}ds\right]+ C_q\mathbb{E}\left[\int_t^{t+h}{( 1+\lVert X_{s}^{\varepsilon } \rVert ^{2q} )}ds\right]\notag
\\&
\le& C_{q,T}h( 1+\lVert  x  \rVert ^{2q}+\lVert  y  \rVert ^{2q} )\notag.
\end{eqnarray}

Putting (\ref{theta1})-(\ref{theta5}) together, the result follows.  \qed
\section{Proof of Main Result} \label{main}
In this section, we intend to give a complete proof for our main result, i.e. the slow component $X^\varepsilon$ converges strongly  to the solution $\bar{X}$ of the averaged equation, as $\varepsilon \rightarrow 0$.
To this end, we
construct a stopping time, i.e.,  for $\varepsilon \in (0,1)$, $ n  >0$,
\begin{eqnarray}\label{stoppingtime}
	 \tau_ n  
	^\varepsilon =\inf \big\{t\ge 0: | {X}_s^\epsilon |_1+ | \bar{X}_s |_1 \ge n  
	\big\} .
\end{eqnarray}
{ \bf {Proof of Theorem \ref{thm2.3}}}.  
The proof is divided into three steps. 

\textbf{Step 1}. For any $q>2$, we know that
 \[
	\mathbb{E}\Big[\sup\limits_{0\leq t\leq T}\big\|X_t^\varepsilon -\bar{X}_t\big\|^{2q}\Big]\leq  \mathbb{E}\left[\mathop{\text{sup}}_{0\le t\le T}\lVert {X}_{t}-\bar{X}_t\rVert ^{2q}\cdot {\chi}_{\{ T>\tau _{ n  }^{\varepsilon } \}}\right]+\mathbb{E}\left[\sup_{0\leq t\leq {T\land \tau _{ n  
	}^{\varepsilon }}}\|X_{t}^{\varepsilon }-\bar{X}_{t}\|^{2q}\right].
	\]
Using the similar argument as in the proof of  Lemma \ref{lem3.1},  Lemma \ref{lem3.33} and  Lemma \ref{lem3.44},  we also get $$\mathbb{E}\left[\underset{0\le t\le T}{\sup} \lVert \bar{X}_{t}^{\varepsilon } \rVert ^{2q}\right] \le C_{q,T}( 1+\lVert  x  \rVert ^{2q}+\lVert  y  \rVert ^{2q}).$$
By H\"{o}lder's inequality, Chebyshev's inequality and Lemma \ref{lem3.1}, we deduce that
 \begin{eqnarray}\label{chebyshev}
&& \mathbb{E}\left[\mathop{\text{sup}}_{0\le t\le T}\lVert {X}_{t}-\bar{X}_t\rVert ^{2q}\cdot {\chi}_{\{ T>\tau _{ n  }^{\varepsilon } \}}\right]\notag\\
&\le& \left[ \mathbb{E}\left(\mathop{\text{sup}}_{0\le t\le T}\lVert {X}_{t}^{\varepsilon }-\bar{X}_t \rVert ^{4q}\right)\right]^{\frac{1}{2}}\big[ \mathbb{P}( T>\tau _{ n  
  }^{\varepsilon } ) \big] ^{\frac{1}{2}}\\
& \le&\frac {C_{q,T}}{n}\left [\mathbb{E}\left(\mathop{\text{sup}}_{0\le t\le T}\lVert {X}_{t}^{\varepsilon } \rVert ^{4q}+\mathop{\text{sup}}_{0\le t\le T}\lVert \bar{X}_t \rVert ^{4q}\right)\right]^{\frac{1}{2}}  \left[ \mathbb{E}\left(\mathop{\text{sup}}_{0\le t\le T}| {X}_{t}^{\varepsilon } |_1+ \mathop{\text{sup}}_{0\le t\le T}| \bar{X}_{s}^{\varepsilon } |_1\right)
\right]^{\frac{1}{2}}\notag\\
&\le& \frac{C_{q,T}}{ n }( 1+\|  x  \|^{2q}+\lVert  y  \rVert ^{2q} )\notag ,
\end{eqnarray}
let $n$ large enough, it yields
 \begin{eqnarray}\label{xjiexbar0}
\lim_{\varepsilon \rightarrow 0} \mathbb{E}\left[\mathop{\text{sup}}_{0\le t\le T}\lVert {X}_{t}-\bar{X}_t\rVert ^{2q}\cdot {\chi}_{\{ T>\tau _{ n  }^{\varepsilon } \}}\right]=0 .
\end{eqnarray}

\textbf{Step 2}.
In this part, we will prove the following result: 
		$$\lim_{\varepsilon \rightarrow 0}\mathbb{E}\left[\sup_{0\leq t\leq {T\land \tau _{ n  
	}^{\varepsilon }}}\|X_{t}^{\varepsilon }-\bar{X}_{t}\|^{2q}\right]=0, ~q>2.$$

According to the definitions of $X_t^\varepsilon$ and $\bar{X}_t$, for any $t\in[0,T]$, we have
\begin{eqnarray}\label{keyxbarxave}
	\mathbb{E}\left[\sup_{0\leq t\leq {T\land \tau _{ n  
	}^{\varepsilon }}}\|X_{t}^{\varepsilon }-\bar{X}_{t}\|^{2q}\right]&\leq& \sum_{i=1}^{3}{I_i}( t ), 
	\end{eqnarray}	
where
\begin{eqnarray*}
	{I_1}( t )&:=&C_q\mathbb{E}\left[\sup_{0\leq t\leq {T\land \tau _{ n  
	}^{\varepsilon }}}\big\|\int_0^t{e}^{( t-s) A}\big(B( X_{s}^{\varepsilon } )- B( \bar{X}_{s} )\big )ds\big\|^{2q}\right],
	\\	
	{I_2}( t )&:=&C_q\mathbb{E}\left[\sup_{0\leq t\leq {T\land \tau _{ n  
	}^{\varepsilon }}}\big\|\int_0^t{e}^{( t-s) A}\big(f_1( X_{s}^{\varepsilon },Y_{s}^{\varepsilon } ) -\bar{f}_1( \bar{X}_{s} )\big)ds\big\|^{2q}\right],\\
	{I_3}( t )&:=&C_q\mathbb{E}\left[\sup_{0\leq t\leq {T\land \tau _{ n  
	}^{\varepsilon }}}\big\|\int_0^t{\int_{|z|<1}{e^{( t-s) A} \big(h_1(X_{s}^{\varepsilon },z)-h_1(\bar{X}_{s},z)\big)}}\tilde{N}_1( ds,dz)\big\|^{2q}\right].
\end{eqnarray*}	
Firstly, estimate the term ${I_1}( t )$,  according to (\ref{qianru}) and Sobolev inequalities in ({\cite{temam1983navier}, Lemma 2.1 ), we get
\begin{eqnarray}\label{J1}
	&&\mathbb{E}\left[\sup_{0\le t\le T\land \tau _{ n  
	}^{\varepsilon }}\lVert I_1( t ) \rVert ^{2q}\right]\notag\\
	&\le&C_q\mathbb{E}\left[ \left(\sup_{0\le t\le T\land \tau _{ n  
	}^{\varepsilon }}\int_0^t{\big( 1+( t-r) ^{-\frac{1}{2}} \big) \big| B( {X}_{r}^{\varepsilon } ) -B( \bar{X}_r) \big|_{-1}dr} \right) ^{2q}\right]\notag\\
&\leq&  C_q\mathbb{E}\left[ \left(\mathop{\text{sup}}_{0\le t\le T\land \tau _{ n  
	}^{\varepsilon }}\int_0^t{\big( 1+( t-r ) ^{-\frac{1}{2}} \big)}\lVert {X}_{r}^{\varepsilon }-\bar{X}_r \rVert \big( | {X}_{r}^{\varepsilon } |_1+| \bar{X}_r |_1 \big) dr\right) ^{2q}\right]\notag
\\
&\leq& C_{q,n}\mathbb{E}\left[ \left(\sup_{0\le t\le T}\int_0^{t\wedge \tau_n^\varepsilon}{( 1+( {t\wedge \tau_n^\varepsilon}-r ) ^{-\frac{1}{2}} )}\lVert {X}_{r}^{\varepsilon }-\bar{X}_r \rVert dr \right) ^{2q}\right]\\
&\leq&C_{q,n}\mathbb{E} \left\{ \sup_{0\leq t\leq T} \left[ \left( \int^{t\wedge \tau^{\varepsilon }_{n} }_{0} \left( 1+\left( t\wedge \tau^{\varepsilon }_{n} -r\right)^{-\frac{1}{2} \cdot \frac{2q}{2q-1} }  \right)  dr\right)^{2q-1}  \int^{t\wedge \tau^{\varepsilon }_{n} }_{0} \| X^{\varepsilon }_{r}-\bar{X}_r \|^{2q} dr\right]  \right\}  \notag
\\&\leq& C_{q,n,T}\int_0^{T}\left[\mathbb{E}\sup_{0\le r\le s\land \tau _{ n  
		}^{\varepsilon }}\lVert {X}_{r}^{\varepsilon }-\bar{X}_r \rVert ^{2q}\right]ds\notag.
\end{eqnarray}
Secondly, estimate the term ${I_2}( t )$, we know that   
\begin{eqnarray*}
{I_2}( t ) &\leq&
C_q\mathbb{E}\left[\sup_{0\leq t\leq {T\land \tau _{ n  
	}^{\varepsilon }}}\big\|\int_0^t{e}^{( t-s) A}\big(f_1( X_{s}^{\varepsilon },Y_{s}^{\varepsilon } ) -f_1(  X_{s(\delta)}^\varepsilon, \hat{Y}_{s}^{\varepsilon })\big)ds\big\|^{2q}\right]\\
&&+
C_q\mathbb{E}\left[\sup_{0\leq t\leq {T\land \tau _{ n  
	}^{\varepsilon }}}\big\|\int_0^t{e}^{( t-s) A}\big(\bar{f}_1( X_{s}^\varepsilon ) -\bar{f}_1( \bar{X}_{s} )\big)ds\big\|^{2q}\right]\\
	&&+
C_q\mathbb{E}\left[\sup_{0\leq t\leq {T\land \tau _{ n  
	}^{\varepsilon }}}\big\|\int_0^t{e}^{( t-s) A}\big(f_1(  X_{s(\delta)}^\varepsilon, \hat{Y}_{s}^{\varepsilon })-\bar{f}_1( X_{s}^\varepsilon)\big)ds\big\|^{2q}\right]\\
	&=:&\sum_{j=1}^{3}{K_j}( t ),
\end{eqnarray*}	
where $s(\delta)=\lfloor\frac{s}{\delta}\rfloor \delta$ is the nearest breakpoint precedings and $\lfloor s\rfloor$  denotes the largest integer which is no more than $s$, and $\hat{Y}^\varepsilon$ is an auxiliary process with the initial value $\hat{Y}_0^\varepsilon = y ^\varepsilon = y $.  For any $t\in\big[k\delta,\min\{(k+1)\delta,T\}\big),$ $k\in \mathbb{N},$ 
	\begin{eqnarray}\label{ytjie}
	\hat{Y}_t^\varepsilon &=&Y_{k\delta}^\varepsilon +\frac{c}{\varepsilon }\int_{k\delta}^tA\hat{Y}_s^\varepsilon  ds+\frac{1}{\varepsilon }\int_{k\delta}^tf_2(X_{k\delta}^\varepsilon ,\hat{Y}_s^\varepsilon )ds+\frac{1}{\sqrt{\varepsilon }}\int_{k\delta}^tdW_s^{Q_2}\notag\\
	&&+\int_{k\delta}^t\int_{{|z|<1}}h_2(X_{k\delta}^\varepsilon ,\hat{Y}_s^\varepsilon, z)\tilde{N}_2^\varepsilon (ds,dz).
	\end{eqnarray}
Applying It\^o's formula,  there exist $\eta 
>0$, such that
\begin{eqnarray*}
	&&\frac{d}{dt}\mathbb{E}\left[\|Y_{t}^{\varepsilon }-\hat{Y}_{t}^{\varepsilon }\|^{2q} \right]
	\\&&=\frac{2qc\lambda _1}{\varepsilon }\mathbb{E}\left[\lVert Y_{t}^{\varepsilon }-\hat{Y}_{t}^{\varepsilon }\rVert ^{2q-2}\big( -\| Y_{t}^{\varepsilon }-\hat{Y}_{t}^{\varepsilon } \|^{2} \big)\right]\\
	&&\quad+\frac{2q}{\varepsilon }\mathbb{E}\left[\lVert Y_{t}^{\varepsilon }-\hat{Y}_{t}^{\varepsilon }\rVert ^{2q-2}\big\langle f_2( X_{t}^{\varepsilon },Y_{t}^{\varepsilon } ) -f_2( X_{k\delta}^{\varepsilon },\hat{Y}_{t}^{\varepsilon } ) ,Y_{t}^{\varepsilon }-\hat{Y}_{t}^{\varepsilon } \big\rangle\right]\\
	&&\quad+\frac{1}{\varepsilon }\mathbb{E}\left[\int_{|z|<1}{\Big( \big\lVert Y_{t}^{\varepsilon }-\hat{Y}_{t}^{\varepsilon }+\big( h_2( X_{t}^{\varepsilon },Y_{t}^{\varepsilon },z) -h_2( X_{k\delta}^{\varepsilon },\hat{Y}_{t}^{\varepsilon },z) \big) \big\rVert ^{2q}-\big\lVert Y_{t}^{\varepsilon }-\hat{Y}_{t}^{\varepsilon} \big\rVert ^{2q} \Big)}\mu_2( dz)\right]\\
	&&\quad-\frac{2q}{\varepsilon }\mathbb{E}\left[\int_{|z|<1}{\big\lVert Y_{t}^{\varepsilon }-\hat{Y}_{t}^{\varepsilon }\big\rVert ^{2q-2}\big\langle h_2( X_{t}^{\varepsilon },Y_{t}^{\varepsilon },z) -h_2( X_{k\delta}^{\varepsilon },\hat{Y}_{s}^{\varepsilon }, z) ,Y_{t}^{\varepsilon }-\hat{Y}_{t}^{\varepsilon } \big\rangle}\mu_2( dz)\right]\\
	&&\le -\frac{q\eta
	}{\varepsilon }\mathbb{E}\left[\lVert Y_{t}^{\varepsilon }-\hat{Y}_{t}^{\varepsilon }\rVert ^{2q}\right]+\frac{C_{q,T,\alpha}}{\varepsilon}\mathbb{E}\left[\lVert X_{t}^{\varepsilon }-{X}_{k\delta}^{\varepsilon }\rVert ^{2q}\right].
\end{eqnarray*}
The Gronwall's inequality \cite{givon2007strong} and Lemma \ref{lem3.44} yields that
\begin{eqnarray}\label{lem3.6}
	\mathbb{E}\left[\lVert Y_{t}^{\varepsilon }-\hat{Y}_{t}^{\varepsilon }\rVert ^{2q}	\right]
	&\le &\frac{C_{q,T,\alpha}}{q \eta}\delta^{\alpha q}( 1+| x |_\alpha^{2q}+\lVert  y  \rVert ^{2q} ) (1- e^{-\frac{q\eta \delta}{\varepsilon}})\notag\\
	&\leq &C_{q,T,\alpha} \delta^{\alpha q}( 1+| x |_\alpha^{2q}+\lVert  y  \rVert ^{2q} ) .
\end{eqnarray}
For $K_1(t)$,  using the contractive property of semigroup $e^{tA},$  (A2), Lemma \ref{lem3.44} and (\ref{lem3.6}), gives
\begin{eqnarray}\label{K1}
K_1(t) &=&
C_q\mathbb{E}\left[\sup_{0\leq t\leq {T\land \tau _{ n  
	}^{\varepsilon }}}\big\|\int_0^t{e}^{( t-s) A}\big(f_1( X_{s}^{\varepsilon },Y_{s}^{\varepsilon } ) -f_1(  X_{s(\delta)}^\varepsilon, \hat{Y}_{s}^{\varepsilon })\big)ds\big\|^{2q}\right]\notag\\
	&\leq&C_q\mathbb{E}\left[\int_0^{T\land \tau _{ n  
	}^{\varepsilon }}\big\|\big(f_1( X_{s}^{\varepsilon },Y_{s}^{\varepsilon } ) -f_1(  X_{s(\delta)}^\varepsilon, \hat{Y}_{s}^{\varepsilon })\big)\big\|^{2q}ds\right]\notag\\
	&\leq& C_q\mathbb{E}\left[\int_0^{T\land \tau _{ n  
	}^{\varepsilon }}\big(\| X_{s}^{\varepsilon }-X_{s(\delta)}^\varepsilon\|^{2q}+\|Y_{s}^{\varepsilon }- \hat{Y}_{s}^{\varepsilon }\|^{2q}\big)ds\right]\\
	&\leq& C_q\int_0^T\mathbb{E}\left[\big(\| X_{s}^{\varepsilon }-X_{s(\delta)}^\varepsilon\|^{2q}+\|Y_{s}^{\varepsilon }- \hat{Y}_{s}^{\varepsilon }\|^{2q}\big)
	\right]dr\notag\\
	&\leq&C_{q,T,\alpha, n}\delta^{\alpha q}
( 1+| x |_\alpha^{2q}+\lVert  y  \rVert ^{2q} )\notag.
	\end{eqnarray}	
For $K_2(t)$,  thanks to the Lipschitz continuity of $\bar{f}_1(\cdot ),$ we have
\begin{eqnarray}\label{K2}
K_2(t)&=&C_q\mathbb{E}\left[\sup_{0\leq t\leq {T\land \tau _{ n  
	}^{\varepsilon }} }\left\|\int_0^t {e}^{( t-s) A}\left(\bar{f}_1( X_{s}^\varepsilon ) -\bar{f}_1( \bar{X}_{s} )\right)ds\right\|^{2q}\right]\notag\\
	&\leq&\mathbb{E} \left[ \sup_{0\leq t\leq {T}  } \left\| \int^{t\wedge \tau^{\varepsilon }_{n} }_{0} {e}^{(t\wedge \tau^{\varepsilon }_{n} -s)A}  \left( \bar{f}_{1} (X^{\varepsilon }_{s})-\bar{f}_{1} (\bar{X}_{s} )\right)  ds\right\|^{2q}  \right]  \\
&\leq&  C_q\int_0^{T}\mathbb{E}\left[\sup_{0\leq r\leq s\land \tau _{ n  
		}^{\varepsilon } }\rVert X_{r}^{\varepsilon } -\bar{X}_{r}^{\varepsilon }\rVert ^{2q} \right]ds\notag.
\end{eqnarray}	
For $K_3(t)$, set $\tilde{n}_t=\lfloor \frac{t}{\delta} \rfloor$, we write
\begin{eqnarray*}
	K_3(t)&=&C_q\mathbb{E}\left[\sup_{0\leq t\leq {T\land \tau _{ n  
	}^{\varepsilon }}}\big\|\int_0^t{e}^{( t-s) A}\big(f_1(  X_{s(\delta)}^\varepsilon, \hat{Y}_{s}^{\varepsilon })-\bar{f}_1( X_{s}^\varepsilon)\big)ds\big\|^{2q}\right]\\
	&\leq&\mathbb{E}\left[ \sup_{0\leq t\leq {T\land \tau _{ n  
	}^{\varepsilon }}} \left(\|K_3^1(t)\|^{2q}+\|K_3^2(t)\|^{2q}+\|K_3^3(t)\|^{2q}\right)\right],
\end{eqnarray*}
where
\begin{eqnarray*}
	K_3^1(t)&=&\sum_{k=0}^{\tilde{n}_t-1}\int_{k\delta}^{(k+1)\delta}e^{(t-s)A}\big[f_1(X_{k\delta}^\varepsilon ,\hat {Y}_s^\varepsilon )-\bar f_1(X_{k\delta}^\varepsilon )\big]ds,\\
	K_3^2(t)&=&\sum_{k=0}^{\tilde{n}_t-1}{\int_{k\delta}^{( k+1 ) \delta}{e^{( t-s ) A}}}\big[ \bar{f}_1( X_{k\delta}^{\varepsilon } ) -\bar{f}_1( X_{s}^{\varepsilon } ) \big] ds,\\
	K_3^3(t)&=&\int_{\tilde{n}_t\delta}^t{e^{( t-s ) A}}\big[{f}_1\big( X_{\tilde{n}_t\delta}^{\varepsilon },\hat{Y}_{s}^{\varepsilon } \big) -\bar{f}_1( X_{s}^{\varepsilon } )\big] ds.
\end{eqnarray*}
For $K_3^2(t),$ according to Lemma \ref{lem3.44}, we have
\begin{eqnarray}\label{j42}
\mathbb{E}\left[\mathop{\text{sup}}_{0\le t\le {T\land \tau _{ n  
	}^{\varepsilon }}
}\lVert K_{3}^{2}( t ) \rVert ^{2q}\right]&\le &C_q\mathbb{E}\int_0^T
\left[\lVert X_{s}^{\varepsilon }-X_{s( \delta )}^{\varepsilon } \rVert ^{2q} \right]ds\notag\\
&\le&C_{q,T,\alpha,n}\delta^{\alpha q}( 1+| x |_\alpha^{2q}+\lVert  y  \rVert ^{2q} ).
\end{eqnarray}
For $K_3^3(t),$ by the construction of  $\hat{Y}_t^\varepsilon$ and similar arguments as in Lemma \ref{lem3.1}, it is easy to obtain $\sup_{0\leq t\leq T}\mathbb{E}\|\hat{Y}_{t}^{\varepsilon }\|^{2q}\leq C_{q,T}( 1+\lVert  x \lVert^{2q}+\lVert  y  \rVert ^{2q} )$, and then follows from Lemma \ref{lem3.1}, it yields
\begin{eqnarray}\label{j43}
&& \mathbb{E} \Big[\mathop{\text{sup}}_{0\le t\le {T\land \tau _{ n }^{\varepsilon }}}\lVert K_{3}^{3}( t ) \rVert ^{2q}\Big]\notag\\
&\le& \mathbb{E}\Big[\mathop{\text{sup}}_{0\le t\le {T}}\int_{\lfloor \frac{{t\wedge \tau_n^\varepsilon}}{\delta} \rfloor\delta}^{t\wedge \tau_n^\varepsilon}\left\lVert  {e^{( {t\wedge \tau_n^\varepsilon}-s ) A}}\left({f}_1\big( X_{\lfloor \frac{{t\wedge \tau_n^\varepsilon}}{\delta} \rfloor\delta}^{\varepsilon },\hat{Y}_{s}^{\varepsilon } \big) -\bar{f}_1( X_{s}^{\varepsilon } )\right)\right \rVert ^{2q}ds\Big]\\
&\le& \delta \mathbb{E} \left[ \sup_{0\leq t\leq {T\wedge \tau_n^\varepsilon} }  \left( 1+\| X^{\varepsilon }_{\lfloor \frac{t}{\delta} \rfloor\delta}\|^{2q} +\| \hat{Y}^{\varepsilon }_{t} \|^{2q} +\| X^{\varepsilon }_{t}\|^{2q} \right) \right] \notag\\
&\le&  C_{q,T,\alpha,n}\delta ( 1+\lVert  x  \rVert ^{2q}+\lVert  y  \rVert ^{2q} )\notag.
\end{eqnarray}
In order to estimate $\mathbb{E}\left[\mathop{\text{sup}}_{0\le t\le T\wedge \tau_n^\varepsilon}\lVert K_3^1(t)\rVert^{2q}\right]$, we get
\begin{eqnarray}\label{k312q}
\mathbb{E}\left[\mathop{\text{sup}}_{0\le t\le {T\land \tau _{ n  
	}^{\varepsilon }}
}\lVert K_{3}^{1}( t ) \rVert ^{2q}\right]\leq \left\{\mathbb{E}\left[\mathop{\text{sup}}_{0\le t\le {T\land \tau _{ n  
	}^{\varepsilon }}
}\lVert K_{3}^{1}( t ) \rVert ^{2(2q-1)}\right]\right\}^{\frac{1}{2}}\left\{\mathbb{E}\left[\mathop{\text{sup}}_{0\le t\le {T\land \tau _{ n  
	}^{\varepsilon }}
}\lVert K_{3}^{1}( t ) \rVert ^{2}\right]\right\}^{\frac{1}{2}}.
\end{eqnarray}
Thanks to Lemma \ref{lem3.1} and the property of $\hat{Y}_t^\varepsilon$, we have
\begin{eqnarray}\label{k312p1}
&&\mathbb{E}\left[\mathop{\text{sup}}_{0\le t\le {T\land \tau _{ n  
	}^{\varepsilon }}
}\lVert K_{3}^{1}( t ) \rVert ^{2(2q-1)}\right]\notag\\
&\leq& \mathbb{E} \left[ \int^{T\wedge \tau^{\varepsilon }_{n} }_{0} \left( \| f_{1}\left( X^{\varepsilon }_{s\left( \delta \right)  },\hat{Y}^{\varepsilon }_{s}\right)  \|^{2(2q-1)} +\| \bar{f}_{1} \left( X^{\varepsilon }_{s\left( \delta \right)  }\right)  \|^{2(2q-1)} \right)  ds\right] \notag \\
&\leq&C_{q,T} \left\{1+\mathbb{E}\left[\sup_{0\leq t\leq T } \| X^{\varepsilon }_{t}\|^{2(2q-1)} \right]+\mathbb{E}\left[\sup_{0\leq t\leq T} \| \hat{Y}^{\varepsilon }_{t} \|^{2(2q-1)} \right]\right\} \\
&\leq& C_{q,T}(1+\| x \|^{4q}+\| y \|^{4q}).\notag
\end{eqnarray}
Moreover, we also have
\begin{eqnarray}\label{j411} 
&&\mathbb{E}\left[\mathop{\text{sup}}_{0\le t\le T\wedge\tau_n^\varepsilon}\left\lVert K_{3}^{1}\left( t \right)\right \rVert ^{2}\right]\notag\\
&=&\mathbb{E}\left[\mathop{\text{sup}}_{0\le t\le T\wedge\tau_n^\varepsilon}\left\lVert \sum_{k=0}^{\tilde{n}_t-1}{e^{\left( t-\left( k+1 \right) \delta \right) A}\int_{k\delta}^{\left( k+1 \right) \delta}{e^{\left( \left( k+1 \right) \delta -s \right) A}\left[ f_1\left( X_{k\delta}^{\varepsilon},\hat{Y}_{s}^{\varepsilon} \right) -\bar{f}_1\left( X_{k\delta}^{\varepsilon} \right) \right] ds}} \right\rVert ^{2}\right]\notag\\
&\le& C {\lfloor \frac{T}{\delta} \rfloor }\mathbb{E} \left[\sum_{k=0}^{ {\lfloor \frac{T\wedge \tau_n^\varepsilon}{\delta} \rfloor }-1}{\left\lVert \int_{k\delta}^{\left( k+1 \right) \delta}{e^{\left( \left( k+1 \right) \delta -s \right) A}\left[ f_1\left( X_{k\delta}^{\varepsilon},\hat{Y}_{s}^{\varepsilon} \right) -\bar{f}_1\left( X_{k\delta}^{\varepsilon} \right) \right] ds}\right \rVert ^{2}\,\,}\right]
\\
&\le&  C_T\left( \frac{T}{\delta} \right) ^{2}\mathbb{E}\left[  \mathop{\max}_{0\le k\le  \lfloor \frac{T}{\delta} \rfloor -1}\left\lVert \int_{k\delta}^{\left( k+1 \right) \delta}{e^{\left( \left( k+1 \right) \delta -s \right) A}\left[ f_1\left( X_{k\delta}^{\varepsilon},\hat{Y}_{s}^{\varepsilon} \right) -\bar{f}_1\left( X_{k\delta}^{\varepsilon} \right) \right] ds}\right \rVert ^{2}\right]
\notag
\\
&\le&  C_T\left( \frac{\varepsilon}{\delta} \right) ^{2}\mathbb{E} \left[\mathop{\max}_{0\le k\le  \lfloor \frac{T}{\delta} \rfloor -1}\left\lVert \int_0^{\frac{\delta}{\varepsilon}}{e^{\left( \delta -s\varepsilon \right) A}\left[ f_1\left( X_{k\delta}^{\varepsilon},\hat{Y}_{s\varepsilon +k\delta}^{\varepsilon} \right) -\bar{f}_1\left( X_{k\delta}^{\varepsilon} \right) \right] ds}\right \rVert ^{2}\right]\notag\\
&\le& C_T\left( \frac{\varepsilon}{\delta} \right) ^{2}\mathop{\max}_{0\le k\le  \lfloor \frac{T}{\delta} \rfloor -1}\int_0^{\frac{\delta}{\varepsilon}}\int_{r}^{\frac{\delta}{\varepsilon}}\Phi_{k}(s,r)dsdr\notag,
\end{eqnarray}
where
\begin{eqnarray*}
&&\Phi_{k}(s,r)\\
&=&\mathbb{E}\left \langle e^{\left( \delta -s\varepsilon \right)  A}\left( f_{1}\left( X^{\varepsilon }_{k\delta },\tilde{Y}^{\varepsilon }_{s\varepsilon +k\delta } \right)  -\bar{f}_{1} \left( X^{\varepsilon }_{k\delta }\right)  \right)  ,e^{\left( \delta -r\varepsilon \right)  A}\left( f_{1}\left( X^{\varepsilon }_{k\delta },\tilde{Y}^{\varepsilon }_{r\varepsilon +k\delta } \right)  -\bar{f}_{1} \left( X^{\varepsilon }_{k\delta }\right)  \right) \right \rangle \\
&=&\mathbb{E} \left\langle e^{\left( \delta -s\varepsilon \right)  A}\left( f_{1}\left( X^{\varepsilon }_{k\delta },Y^{X^{\varepsilon }_{k\delta },Y^{\varepsilon }_{k\delta }}_{s}\right)  -\bar{f}_{1} \left( X^{\varepsilon }_{k\delta }\right)  \right)  ,e^{\left( \delta -r\varepsilon \right)  A}\left( f_{1}\left( X^{\varepsilon }_{k\delta },Y^{X^{\varepsilon }_{k\delta },Y^{\varepsilon }_{k\delta }}_{r}\right)  -\bar{f}_{1} \left( X^{\varepsilon }_{k\delta }\right)  \right)  \right\rangle ,
\end{eqnarray*}
here the distribution of $\left( X_{k\delta}^{\varepsilon},\hat{Y}_{s+k\delta}^{\varepsilon} \right) $ coincides with  the distribution of $ \left( X_{k\delta}^{\varepsilon},Y_{{s}/{\varepsilon}}^{X_{k\delta}^{\varepsilon},Y_{k\delta}^{\varepsilon}} \right)$. Similar as the argument in \cite{xu2015lp, 2011Strong}, using Lemma \ref{lem3.1}, one can verify that
\begin{eqnarray}\label{gksr}
\Phi_{k} (s,r)\leq C\mathbb{E} \left( 1+\| X^{\varepsilon }_{k\delta }\|^{2} +\| Y^{\varepsilon }_{k\delta }\|^{2} \right)  e^{-\frac{1}{2} \left( s-r\right)  \eta }\leq C_{T}\left( 1+\| x\|^{2} +\| y\|^{2} \right)  e^{-\frac{1}{2} \left( s-r\right)  \eta }.
\end{eqnarray}
Then, thanks to (\ref{j411}) and (\ref{gksr}), we get that for any $\varepsilon\in(0,1)$,
\begin{eqnarray}\label{k312}
\mathbb{E}\left[\mathop{\text{sup}}_{0\le t\le {T\land \tau _{ n  
	}^{\varepsilon }}
}\lVert K_{3}^{1}( t ) \rVert ^{2}\right]
&\le &C_q( \frac{\varepsilon }{\delta} ) ^{2}
\frac{ \delta}{\varepsilon } ( 1+\lVert  x  \rVert ^{2q}+\lVert  y  \rVert ^{2q} )\notag\\
&=&C_{q,T,\alpha,n}\frac{\varepsilon }{\delta}\big( 1+\lVert  x  \rVert ^{2q}+\lVert  y  \rVert ^{2q} \big).
\end{eqnarray}
Substituting (\ref{k312}) and (\ref{k312p1}) into (\ref{k312q}), we get
\begin{eqnarray}\label{j41}
\mathbb{E}\left[\mathop{\text{sup}}_{0\le t\le {T\land \tau _{ n  
	}^{\varepsilon }}
}\lVert K_{3}^{1}( t ) \rVert ^{2q}\right]
&=&C_{q,T,\alpha,n}\sqrt{\frac{\varepsilon }{\delta } } \big( 1+\lVert  x  \rVert ^{2q}+\lVert  y  \rVert ^{2q} \big).
\end{eqnarray}
Combining (\ref{j42}), (\ref{j43}) and (\ref{j41}), we can conclude
\begin{eqnarray}\label{K3}
\mathbb{E}\left[\mathop{\text{sup}}_{0\le t\le {T\land \tau _{ n  
	}^{\varepsilon }}
}\lVert K_3(t)\rVert ^{2q}\right]
\le C_{q,T,\alpha,n}(\delta^{\alpha q}+\delta +\sqrt{\frac{\varepsilon }{\delta } })( 1+| x |_\alpha^{2q}+\lVert  y  \rVert ^{2q} ).
\end{eqnarray}
 Thirdly, estimate the term ${I_3}(t)$, by the Kunita's first inequality, the H\"{o}lder inequality, (A2) and Young's inequality, we have
 \begin{eqnarray}\label{J3}
{I_3}( t ) &\leq&C_q\mathbb{E}\left[\sup_{0\leq t\leq {T\land \tau _{ n  
	}^{\varepsilon }}}\big\|\int_0^t{\int_{|z|<1}{e^{( t-s) A} \big(h_1(X_{s}^{\varepsilon },z)-h_1(\bar{X}_{s},z)\big)}}\tilde{N}_1( ds,dz)\big\|^{2q}\right]\notag\\
&\leq& C_q \mathbb{E}\left[ \int_{0}^{T\land \tau _{ n  
		}^{\varepsilon }}\int_{|z|<1}\|h_1( X_{t}^{\varepsilon },z ) -h_1( \bar{X}_{t},z )  \|^2\mu_{1}( dz ) dt\right] ^q \\
&&+ C_q\mathbb{E}\left[ \int_{0}^{T\land \tau _{ n  
		}^{\varepsilon }}\int_{|z|<1}\| h_1( X_{t}^{\varepsilon },z ) -h_1( \bar{X}_{t},z)  \|^{2q}\mu_{1}( dz ) dt\right] \notag\\\notag
&\leq&  C_q\int_0^{T}\mathbb{E}\left[\sup_{0\leq r\leq s\land \tau _{ n  
		}^{\varepsilon } }\lVert X_{r}^{\varepsilon } -\bar{X}_{r}\rVert ^{2q}\right]ds\notag.
\end{eqnarray}
Again, substituting (\ref{J1}), (\ref{K1}),  (\ref{K2}), (\ref{K3}) and (\ref{J3}) into (\ref{keyxbarxave}), we can derive
\begin{eqnarray*}
	\mathbb{E}\left[\mathop{\text{sup}}_{0\le t\le T\land \tau _{ n  
	}^{\varepsilon }}\lVert {X}_{t}^{\varepsilon }-\bar{X}_t \rVert^{2q}\right]
&\le& C_{q,T,\alpha,n}(\delta^{\alpha q}+\delta +\sqrt{\frac{\varepsilon }{\delta } } )( 1+| x |_\alpha^{2q}+\lVert  y  \rVert ^{2q} )\\
&&+C_{q,T,\alpha,n} \int_0^T\mathbb{E}\left[\mathop{\text{sup}}_{0\le r\le s\land \tau _{ n  
		}^{\varepsilon }}\lVert {X}_{r}^{\varepsilon }-\bar{X}_r \rVert ^{2q}\right]ds.
\end{eqnarray*}
Using Gronwall's inequality, we get
\begin{eqnarray*}
	\mathbb{E}\Big[\mathop{\text{sup}}_{0\le t\le T\land \tau _{ n  
		}^{\varepsilon }}\lVert {X}_{t}^{\varepsilon }-\bar{X}_t \rVert^{2q}\Big] &\le &C_{q,T,\alpha,n}( 1+| x |_\alpha^{2q}+\lVert  y  \rVert ^{2q} ) (\delta^{\alpha q}+\delta+\sqrt{\frac{\varepsilon }{\delta } })e^{C_{q,T,\alpha,n} }.
\end{eqnarray*}
 Taking   $\delta=\varepsilon ^{\frac{1}{2}}$,  we conclude
  \begin{eqnarray}\label{xjiexbartau}
\lim_{\varepsilon \rightarrow 0} \mathbb{E}\left[\underset{_{0\le t\le T}}{\text{sup}}\lVert {X}_{t}^{\varepsilon }-\bar{X}_t \rVert ^{2q}\cdot {\chi}_{\{ T\le \tau _{ n  
  	}^{\varepsilon } \}}\right]=0.
  \end{eqnarray}   
  
\textbf{Step 3}. Taking  $\delta=\varepsilon^{\frac{1}{2}}$ and letting $\varepsilon\rightarrow 0$ first, then $n\rightarrow \infty$, by (\ref{chebyshev}) and  (\ref{xjiexbartau}),  for $q>2$, we have 
\[
		\lim_{\varepsilon  \rightarrow 0}\mathbb{E}\Big[\sup\limits_{0\leq t\leq T}\big\|X_t^\varepsilon -\bar{X}_t\big\|^{2q}\Big]=0.  
		\]
 Applying  the Jensen inequality,  for $2\leq p\leq q$, we get
\begin{eqnarray*}
\mathbb{E}\left[\mathop{\text{sup}}_{0\le t\le T
}\lVert {X}_{t}^{\varepsilon }-\bar{X}_t  \rVert ^{p}\right]\leq \left \{ \mathbb{E}\left[\mathop{\text{sup}}_{0\le t\le T
}\lVert {X}_{t}^{\varepsilon }-\bar{X}_t  \rVert ^{2q}\right] \right \}^{\frac{p}{2q}},
\end{eqnarray*} 
thus, taking  $\delta=\varepsilon^{\frac{1}{2}}$,  we have
 \begin{eqnarray} \label{pge1}
\lim_{\varepsilon \rightarrow 0}\mathbb{E}\left[\mathop{\text{sup}}_{0\le t\le T}\lVert {X}_{t}^{\varepsilon }-\bar{X}_t \rVert ^{p}\right]=0, ~p \ge 2.
\end{eqnarray} }
This completes the proof of Theorem \ref{thm2.3}.   \qed

\section{Numerical Simulations}\label{numerical}
In this section, we present several examples to illustrate graphically the averaging principle for system (\ref{orginal0}). Consider the Burgers equation with an external force
\begin{eqnarray} \label{example0}
\frac{\partial u}{\partial t}=\frac{\partial ^2u}{\partial \xi^2}+u\frac{\partial u}{\partial \xi} +F(u, v, \zeta),
\end{eqnarray}
The external force $F ( u,  v, \zeta )$ depends on $u$, $v$ and $\zeta$, which makes the fluid run complicatedly. As we know, the term $v$ can be governed by a fast evolution process, and $\zeta$ may be a stochastic process. Specifically, taking 
\begin{eqnarray*}
\begin{aligned}
	&F( u,  v, \zeta)dt =-( u+v)dt+\zeta dt,\\
	&dv= -\frac{1}{\varepsilon}v dt+\frac{1}{\sqrt{\varepsilon}}dW_t,\\
	&\zeta dt= dW_t+ \int_{|z|<1}uz\tilde{N}( dt,dz ),
\end{aligned}
\end{eqnarray*}
then the equation (\ref{example0}) can be rewritten as follows with suitable initial data and Dirichlet boundary condition
\begin{eqnarray}\label{example}
\left\{ 
	\begin{array}{l}
		d{u}_{t}^{\varepsilon}( \xi ) =\Big[ \Delta {u}_{t}^{\varepsilon}( \xi ) +\frac{1}{2} \frac{\partial}{\partial \xi}\big( {u}_{t}^{\varepsilon}( \xi ) \big) ^2 -( {u}_{t}^{\varepsilon}( \xi)+ {v}_{t}^{\varepsilon}( \xi)) \Big] dt+dW_t+\int_{|z|<1}{{u}_{t}^{\varepsilon}( \xi )}z\tilde{N}( dt,dz ), \\
d{v}_{t}^{\varepsilon} = -\frac{1}{\varepsilon}{v}_{t}^{\varepsilon}dt+\frac{1}{\sqrt{\varepsilon}}dW_t,\\
{u}_{0}^{\varepsilon} =2, {v}_{0}^{\varepsilon}=1, ~\xi \in [ 0, 1],\\
{u}_{t}^{\varepsilon}( 0 ) ={u}_{t}^{\varepsilon}( 1 ) ={v}_{t}^{\varepsilon}( 0 ) ={v}_{t}^{\varepsilon}\left( 1 \right) =0, t \in [ 0, 1],
	\end{array} \right. 
\end{eqnarray}
which is a typical example of the system (\ref{orginal0}). It is not hard to deduce the averaged equation associated with the system (\ref{example})
\begin{eqnarray}\label{averageequ2}
d{\bar{u}}_{t}(\xi) =\left[ \Delta {\bar{u}}_{t}(\xi) + \frac{1}{2}\frac{\partial}{\partial \xi}\big( {\bar{u}}_{t}( \xi )\big ) ^2-{\bar{u}}_{t}( \xi)\right] dt+dW_t+\int_{|z|<1}{{\bar{u}}_{t}}( \xi )z\tilde {N}_1(dt,dz), \quad {\bar{u}}_{0}(\xi)=2.
\end{eqnarray}
As shown in Fig.\ref{fig}, the solution ${u}_{t}^{\varepsilon}( \xi )$ of system (\ref{example}) converges to the solution ${\bar{u}}_{t}(\xi)$ of the averaged equation (\ref{averageequ2}) if $\varepsilon\rightarrow 0$ in $L^{\frac{p}{2}}(\Omega, L^{2}(0, 1))$, $p=3, 4$, which is in accordance with the conclusion given in Theorem \ref{thm2.3}.

\begin{figure}[htbp]
\centering
\subfigure[$p=3$]{      
\includegraphics[width=7cm]{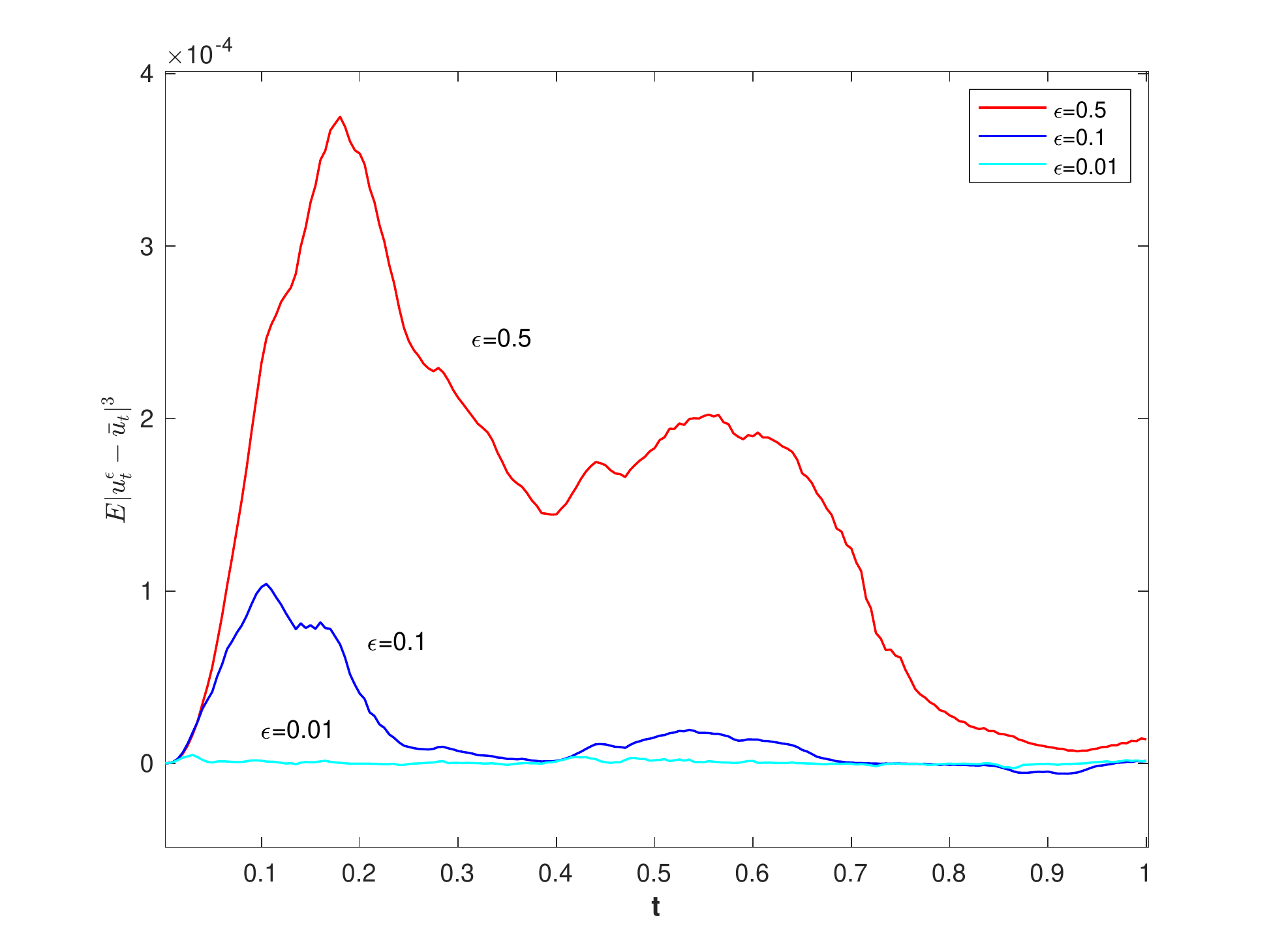}
}
\hspace{0in}
\subfigure[$p=4$]{
\includegraphics[width=7cm]{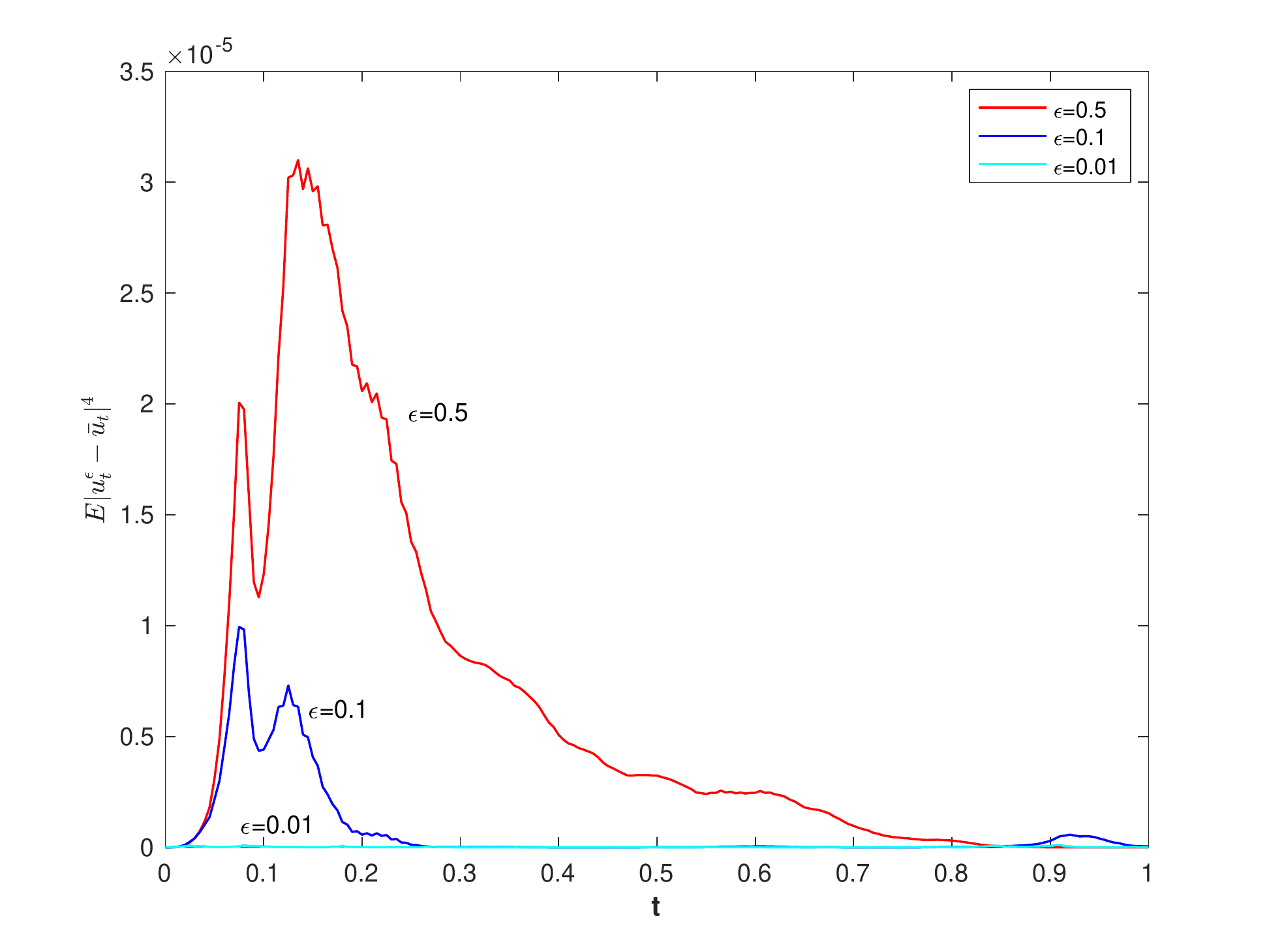}
}
\caption{Convergence for one dimensional stochastic Burgers equation as $\varepsilon$ goes to zero.}
\label{fig}
\end{figure}

\section*{Acknowledgments}
 
 This work was partly supported by the Key International (Regional) Cooperative Research Projects of the NSF of China (Grant 12120101002), the NSF of China (Grant 12072264, 11802236), the Fundamental Research Funds for the Central Universities, the Research Funds for Interdisciplinary Subject of Northwestern Polytechnical University,  the Shaanxi Provincial Key R\&D Program (Grants 2020KW-013, 2019TD-010).


\bibliography{references}

\end{document}